\documentclass[pdflatex,sn-mathphys-num]{sn-jnl}% Math and Physical Sciences Numbered Reference Style
\usepackage{helvet}
\usepackage{url}              %
\usepackage{xcolor}           % 
\usepackage{ulem}             % 
\usepackage{tikz}             %      

\usepackage{pdflscape}        % 
\usepackage{colortbl}         % 
\usepackage{wrapfig}          % 
\usepackage{graphicx}         % 
\usepackage{etoolbox}
\usepackage{todonotes}
\usepackage{amsmath}
\usepackage{amssymb}    
\usepackage{caption}
\usepackage{subcaption}
\usepackage{scalerel}
\usepackage{transparent}
\usepackage{adjustbox}
\usepackage[numbers]{natbib}
\usepackage{doi}
\usetikzlibrary{arrows.meta}

\begin{document}

\title[Article Title]{FSI modeling of case-specific nonlinear carotid artery mechanics and the role of outlet boundary conditions}


\author*[1]{\fnm{Tristan} \sur{Probst}}\email{t.probst@rptu.de}

\author[1]{\fnm{Anna} \sur{Hundertmark}}

\author[2]{\fnm{Ashkan} \sur{Shiravand}}
\author[2]{\fnm{Giorgio} \sur{Cattaneo}}

\affil*[1]{\orgdiv{Institute of Mathematics}, \orgname{University of Kaiserslautern-Landau}, \orgaddress{\street{Fortstr. 7}, \city{Landau}, \postcode{76829}, \state{Rhineland-Palatinate}, \country{Germany}}}

\affil[2]{\orgdiv{Institute of Biomedical Engineering}, \orgname{University of Stuttgart}, \orgaddress{\street{Seidenstr. 36}, \city{Stuttgart}, \postcode{70174}, \state{Baden-Württemberg}, \country{Germany}}}


\abstract{ Compliant arterial wall mechanics, contributing to the Windkessel effect in the carotid artery (CA), have significant impact on hemodynamic patterns and surface shear indicators for cardiovascular diseases, e.g. atherosclerosis. 
    In this study, we extend the linear elastic Fluid-Structure Interaction (FSI) framework to account for nonlinear strain-dependent wall behavior. The model incorporates a Young’s modulus generalized from tensile tests on silicone phantoms to capture the nonlinear stress-strain relation of arterial tissue.
    The final computational model, using a resistance-type boundary condition  and in vitro measured stress-strain relation is validated  against both, in vitro assessed silicon CA phantom as well as published clinical data on the flow splitting to the daughter branches in CA bifurcations. 
    To better reflect physiological conditions, our model is subsequently extended to incorporate prestress of patient-specific geometries, and clinically measured stress-strain relations, followed by validation against clinical CA data.
    The present study demonstrates the feasibility of strain-dependent Young’s (elastic) modulus as a means to enhance the capacity of the linear elastic framework to accurately represent the physiologically nonlinear mechanics of arterial walls, striking a balance between implementation effort and physiological fidelity. This approach yields realistic strains, volumetric inflation behavior and nonlinear pressure-volume relationships. Furthermore, the study reveals the importance of proper outlet boundary conditions and the shortcomings of resistance-type boundary condition in patient-specific modeling, leading to non-physiological pressure profiles and non-physiological temporal flow splitting.}




\maketitle
%-------------------------------------------------------------------------------------------------------------------------
\section{Introduction}
Over the past few decades, numerical modeling of cardiovascular flow has emerged as a prominent and popular topic within the field of computational fluid dynamics.
Disturbed flow patterns and shear stress predictors evaluated from numerical simulations serve as hemodynamic indicators.
High, as well as low, and oscillating wall shear stresses (WSS)
\cite{arzani, Giorgio_low_WSS, Giorgio_WSS2, taylor1998finite, Tremmel_OSI, ku} 
and their multi-directional behavior, measured e.g. by transversal or longitudinal WSS or cross flow index \cite{hoogendorn2020multidirectional, mohamied2015understanding,  pfeifer2013computation, pfeifer2013low, morbiducci2015rational,ProRiHu23,kok2019influence,carpenter2023nonlinear},  have been linked  to atherosclerotic promotion.
The vessel wall possesses the capacity to accumulate kinetic energy and subsequently release it, which regulates and smooths the blood dynamics in the arteries. 
Fluid-structure interaction (FSI) models describe the interplay between fluid and solid tissue dynamic. Several constitutive laws have been applied and discussed in the context of vessel wall mechanics.
In~\cite{colciago2014comparison} 3D elastodynamic and reduced order, thin-wall or membrane FSI models have been compared. More complex and multilayered hyperelastic and anisotropic fiber-reinforced models for vessel tissue modeling can be found in~\cite{holzapfel2000constitutive}.
The comparisons between easily available linearly elastic and hyper-elastic constitutive laws for wall dynamic modeling of a cerebral aneurysm in~\cite{torii2008fluidstructure} indicate some overshooting of deformation in the case of linear-elasticity model.
The displacement patterns for both models, however, have been found to be similar, indicating the usefulness of both material models.
Vo{\ss} et al. \cite{voss2016fluid} also highlight the importance of inhomogeneous wall thickness regarding the WSS surface distribution in the FSI modeling of compliant intracranial aneurysm.
In view of a significant impact of arterial wall deformability on hemodynamic patterns and surface risk indicators, the incorporation of realistic compliance modeling is of paramount importance in the hemodynamic validation process.
In the context of wall compliance and the clinically observed \textit{Windkessel effect},  FSI simulations should accurately replicate the crucial function of the arterial wall as an energy storage reservoir. 

The present study focuses on numerical compliance modeling based on our previous validation  of the compliance behavior with an in vitro physical model \cite{ShiRiHu24} and on clinical stress-strain behavior from the literature \cite{faturechi2019mechanical}. 
The physical model represents a patient-specific carotid artery  with a stenotic occlusion region and is fabricated  using material commonly used in  
in-vitro modelling and validation, 
a silicone -\textit{Dragon Skin 10 NV} (DS), see e.g. \cite{shiravand2025core, shiravand2026}. For fabrication reasons the wall thickness of the artery phantom was artificially but consistently increased and some artery branches has been cropped.
\textit{Diametric distensibility}, \textit{volumetric compliance} as well as  \textit{pressure-volume curves} have been compared with numerical simulation and extend the results of \cite{ShiRiHu24}.
As confirmed in \cite{ShiRiHu24} by tensile tests and in \cite{faturechi2019mechanical} by clinical observation, the realistic tensile behavior of the phantom materials as well as of human artery walls exhibits nonlinear stress-strain relations within physiological pressure ranges.
In order to reproduce this behavior and to achieve better conformance with the in vitro experiments without resorting to more complex e.g., hyperelastic models, we adapt our FSI model based on linear elasticity by introducing a dependency of the elasticity (Young's) modulus on a specific strain metric. 

The paper is organized as follows.
In section \ref{sec:model}, we provide a comprehensive description of the mathematical model including the strain metric employed for the generalization of the Young's modulus. This approach facilitates more precise reconstruction of the compliance behavior, whose validation is summarized in section \ref{sec:sealedcarotis} for sealed artery for better comparability with the in vitro measurements set up. The convergence of the numerical solution towards a highly resolved reference solution, confirming respective orders of convergence as well as the convergence of the artery volume and of the volumetric compliance with respect to mesh resolution are documented in section \ref{sec:convergencestudy}.

In section \ref{sec:openartery} we transfer our FSI modeling to the unsealed and more complex geometry of the patient's carotid artery tree with several branching vessels and realistic wall thickness,  allowing the downstream flow under consideration of resistance-inertance boundary condition (BC).
Based on the Laplace formula introduced in \cite{fung1993biomechanics}, the elasticity module for the thinner wall geometry is %increased 
calibrated by a constant factor necessary to preserve the compliance behavior, i.e., the relation between  volume increase  and  pressure action.
The final realistic unsealed, patient-specific model is validated against in vitro assessed pressure-volume behavior, as in section \ref{sec:sealedcarotis} as well as clinical flow division data  for carotid artery (CA) bifurcations, revealing conforming compliance behavior for our generalized Young's modulus approach.
Concerning the flow division ratios between the CA branches and the pressure behavior, results indicate the necessity of application of an advanced version of the here presented resistance-type outlet BC.

Finally, the strictly linear FSI model is benchmarked against the generalization of Young's modulus in section \ref{subsec:Invivo} for  a physiological stress-strain curve derived from individual CA samples (in vivo data), reported by Faturechi et al. \cite{faturechi2019mechanical}. Here, prestress is incorporated into the model to account for the prestressed configuration of patient-specific arteries extracted through clinical imaging.  This comparison underscore the enhanced capability of the adapted approach to capture different mechanical behaviors, and therefore to support and validate in vitro experiments as well as accurately capturing physiological patient specific arterial compliance behavior.
%-------------------------------------------------------------------------------------------------------------------------
\section{Models and methods}
\label{sec:model}

\begin{figure}[htpb]
    \centering
    \begin{subfigure}[b]{0.48\textwidth}
        \begin{tikzpicture}
        \node[anchor=south west, inner sep=0] (img) at (0,0) {%
            \adjincludegraphics[trim={0 0 {0.4\width} 0}, clip, height=8cm, width=\textwidth]{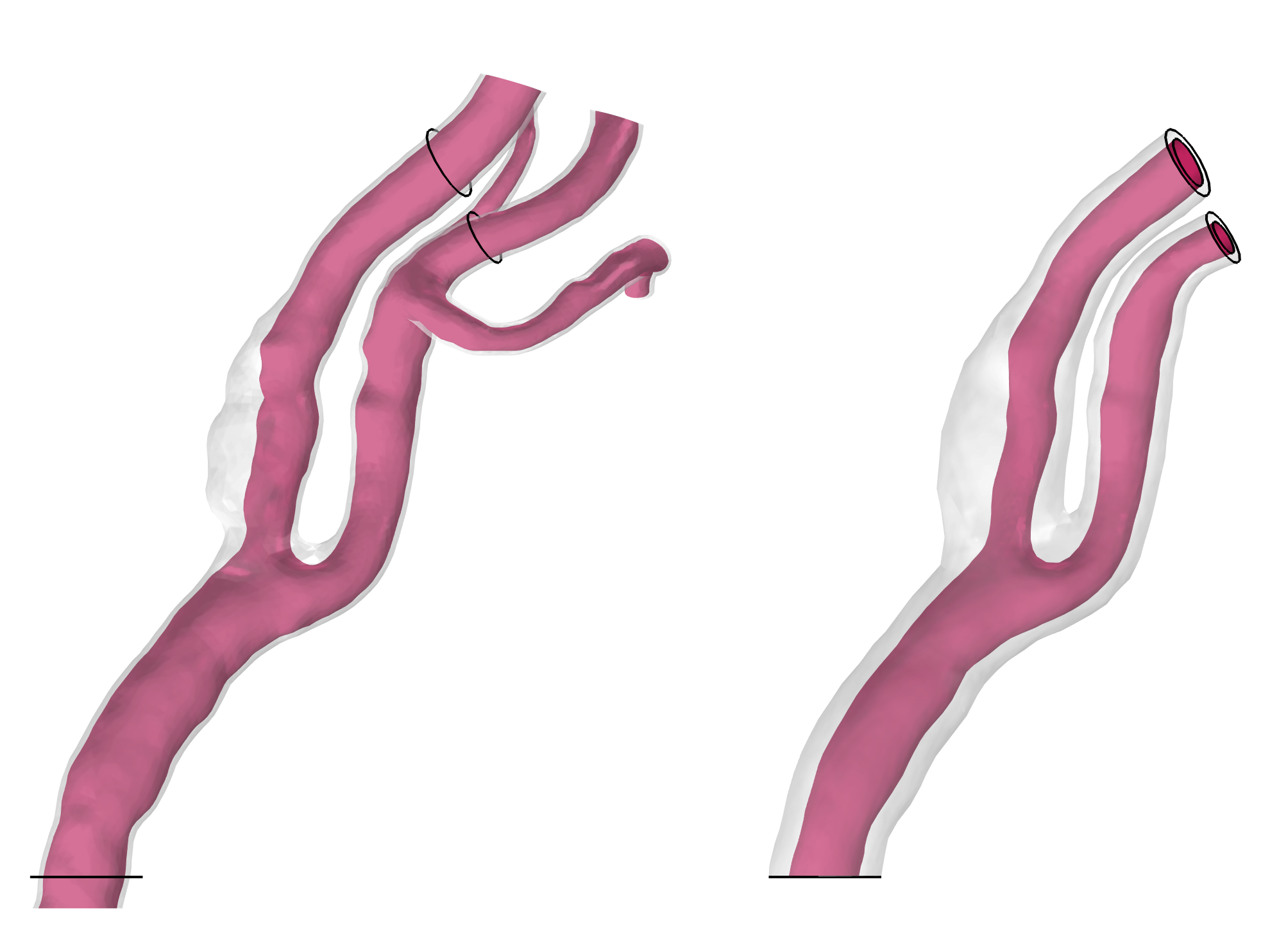}%
        };
        \begin{scope}[x={(img.south east)}, y={(img.north west)}]
            
            \tikzset{
                myarrow/.style={-{Stealth[scale=1.2]}, thick, draw=black},
                labelstyle/.style={ fill=white, fill opacity=0.85, text opacity=1, inner sep=2pt, rounded corners=2pt}
            }
            \tikzset{
                myarrow/.style={-{Stealth[scale=1.2]}, thick, draw=black, rounded corners=3pt},
                labelstyle/.style={ inner sep=2pt}
            }
            %Inflow
            \node[labelstyle, anchor=west] (inflow) at (0.25, 0.03) {$\Gamma_{in}$};
            \draw[myarrow] (0.1, 0.05) --  (0.1, 0.03) --(inflow.west)  ;
            %fsi
            \node[labelstyle, anchor=west] (fsi) at (0.38, 0.25) {$\Gamma_{fsi}$};
            \draw[myarrow]  (0.22, 0.25) -- (fsi.west);
            %outflow
            \node[labelstyle, anchor=west] (outflow) at (0.9, 0.9) {$\Gamma_{out}$};
            \draw[myarrow]  (0.7, 0.9) -- (outflow.west);
            \draw[myarrow] (0.8,0.87) -- (0.8, 0.9) -- (outflow.west) ;
            \draw[myarrow] (0.85,0.7) -- (0.9, 0.7) -- (outflow.west) ;
        \end{scope}
    \end{tikzpicture}
        \caption{}
        \label{fig:realistic_geometry}
    \end{subfigure}
    \hfill
    \begin{subfigure}[b]{0.48\textwidth}
        \centering
             \begin{tikzpicture}
        \node[anchor=south west, inner sep=0] (img) at (0,0) {%
            \adjincludegraphics[trim={{0.6\width} 0 0 0}, clip, height =8cm,width=0.8\textwidth]{RealisticAndSimpleGeometry.png}
            };

        \begin{scope}[x={(img.south east)}, y={(img.north west)}]

            % Style-Definition für Pfeile & Textboxen
            \tikzset{
                myarrow/.style={-{Stealth[scale=1.2]}, thick, draw=black},
                labelstyle/.style={ fill opacity=0,inner sep=2pt, rounded corners=2pt}
            }
            % \draw[help lines, red, step=0.1] (0,0) grid (1,1); \foreach \x in {0,0.2,...,1} \node[red] at (\x,0) {\tiny \x}; \foreach \y in {0,0.2,...,1} \node[red] at (0,\y) {\tiny \y};
            \tikzset{
                myarrow/.style={-{Stealth[scale=1.2]}, thick, draw=black, rounded corners=3pt},
                labelstyle/.style={ fill=white, fill opacity=0.85, text opacity=1, inner sep=2pt}
            }
             \tikzset{
                myarrowt/.style={-{Stealth[scale=0.9]}, thin, draw=black, rounded corners=3pt},
                labelstyle/.style={ fill=white, fill opacity=0.85, text opacity=1, inner sep=2pt}
            }
            %Inflow
            \node[labelstyle, anchor=west] (inflow) at (0.25, 0.03) {$\Gamma_{in}$};
            \draw[myarrow]  (0.1, 0.08) -- (0.1, 0.03) -- (inflow.west);
            %fsi
            \node[labelstyle, anchor=west] (fsi) at (0.73, 0.25) {$\Gamma_{fsi}$};
            \draw[myarrow] (0.25, 0.25) -- (fsi.west);
            \draw[myarrow] (0.82, 0.82) -- (0.82,0.28);
            %outflow
            \node[labelstyle, anchor=north] (outflow) at (0.9, 0.95) {sealed outlets};
            %ECA ICA Beschriftung
            \node[ inner sep=2pt, anchor=north] (ICA) at (0.5, 0.5) {ICA};
            \node[ inner sep=2pt, anchor=north] (ICA) at (0.7, 0.5) {ECA};
        \end{scope}
    \end{tikzpicture}
        \caption{}
        \label{fig:cropped_geometrie}
    \end{subfigure}

    \caption{ a) Original CA geometry with four outlet cross sections and realistic wall thickness. b) Simplified geometry with thick walls and two main branches ICA and ECA, wehre remaining branches have been cropped.}
    \label{fig:geometries}
\end{figure}
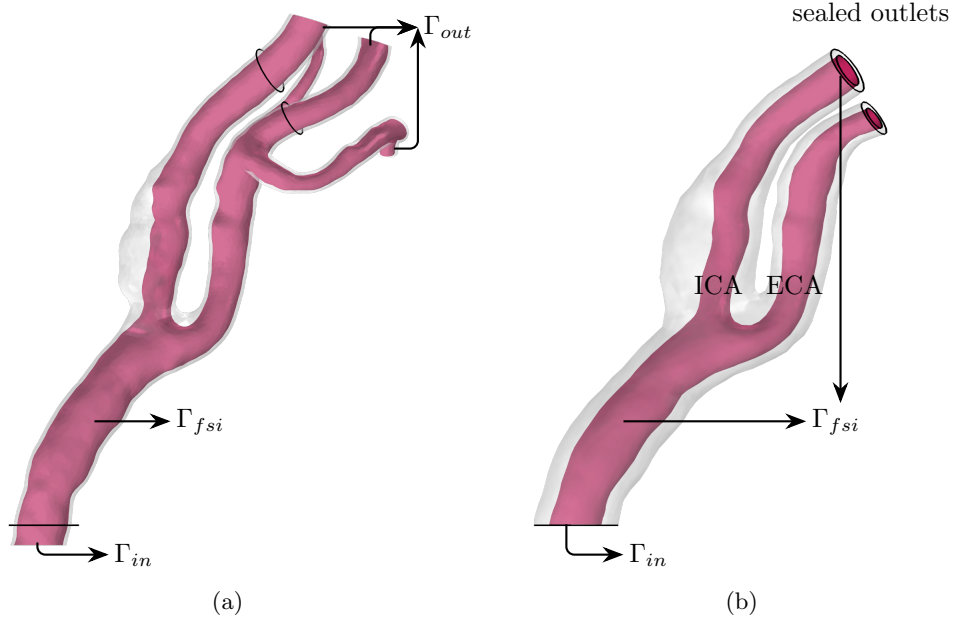

A 3D model of an carotid artery is derived from CTA scans of a single individual patient. 
Atherosclerotic plaque and regions of wall thinning are incorporated leading to varying artery thickness in these areas.
The geometric reconstruction process from anatomical data is described in our previous work \cite{ShiRiHu24}.
In this contribution we consider two different artery geometries: the first is a high-fidelity model of the patient's artery tree structure with accurate sub-branching of the internal (ICA) and external carotid artery (ECA),  see Fig \ref{fig:realistic_geometry}. 
A second, simplified model is used to analyse the compliance behaviour by numerical simulations and compare it with in vitro  experimental studies for the fabricated physical arterial model with sealed outlets.
The second geometry neglects anatomic details for externa branches and consist only of the internal and  the main external vessel branches ICA and ECA, see Fig. \ref{fig:cropped_geometrie}, with  artificially 
thickened vessel wall for consistency with the laboratory compliance study.
Indeed, a sufficient wall thickness decrease the risk of the  occurrence of fissures in the molding process of artery phantom fabrication, which poses a challenge for the laboratory setup.
The outlets of the simplified geometry are closed by
 a thick slice of compliant solid  material to seal the geometry.
 The thickness of the closure outlet caps is comparable to that of the  remaining vessel wall, compare also Fig. \ref{fig:straindistribution} ,resulting in a computational domain with one inlet surface and compliant vessel walls that can inflate longitudinally and transversely. 

\subsection{Numerical model for sealed artery} \label{subsec:sealedsetup}
We denote with $\Omega_t=\Omega^f_t \cup \Omega^s_t\;,\;t\in[0,T]$ the deforming fluid and structural domain at time point $t$, respectively 
 
and formulate the mathematical model of the fluid-structure interaction problem as follows.

\subsubsection*{Fluid flow subproblem}
The fluid flow is modeled by incompressible  Navier-Stokes equations considering Newtonian fluid. 
The Arbitrary Lagrangian-Eulerian (ALE) formulation takes into account the spatial deformation %$x(X,t)$ 
of the reference domain (of mesh-points) with actual and reference coordinates denoted as $x(X,t)\in \Omega_t^f, X\in \Omega_0^f$ respectively. In the mixed ALE formulation the time derivative is expressed in reference, i.e. material (Lagrangian) coordinates ($X$) for the purpose of time-discretization reasons, that is compensated by additional convective term related to the deformation velocity $\boldsymbol{w}$,

\begin{equation}\label{eq:ns}
\begin{gathered}
    \rho_f\frac{D \boldsymbol{u}}{D t}+ \rho_f\big((\boldsymbol{u}-\boldsymbol{w}) \cdot  \nabla\big) \boldsymbol{u} 
    = div \ \boldsymbol{T}_f, \quad div \ \boldsymbol{u}=0 \qquad \text{in}\quad \Omega^f_t.\\
    \boldsymbol{T}_f=-p\boldsymbol{I}+2\mu\big(\nabla \boldsymbol{u}+ \nabla \boldsymbol{u} ^T\big)
    \\
    \boldsymbol{\tilde{w}}(X,t):=\frac{\partial x(X,t)}{\partial t},\  X\in  \Omega^f_0
\end{gathered}
\end{equation}
%}
% 
Here $\boldsymbol{u}$ is the fluid velocity; $\boldsymbol{T}_f$ the Cauchy stress tensor for the fluid including a symmetric strain tensor.  $\frac{D\boldsymbol{u}}{Dt}(x,t)=\partial_t \tilde{\boldsymbol{u}}(X,t)=\frac{\partial \boldsymbol{u}(x,t)}{\partial t}+\nabla\boldsymbol{u}(x,t)\cdot \boldsymbol{w}(x,t)$  is the material time derivative with Lagrangian fluid velocity $\tilde{\boldsymbol{u}}(X,t)$, and $\tilde{\boldsymbol{w}}(X,t)$, describes the  deformation velocity of the reference domain $\Omega^f_0$.\\

The  surface of the fluid domain $\partial \Omega^f_t=:\Gamma^{f}=\Gamma^f_{in}\cup\Gamma^t_{fsi}$ of the sealed artery  consists of the inlet surface and the fluid-structure-interaction interface.
We assume the zero initial velocity and consider a  constant normal inflow velocity $v_0$ at the inlet surface prescribing
\begin{equation}
    \boldsymbol{u}\equiv r(t) v_{0}\cdot\boldsymbol{n}_f, \quad  v_0=\tfrac{|\Omega^f_0|}{T\,|\Gamma^f_{in}|} \quad \text{ on } \Gamma_{in}^f,
\end{equation}
 which produces a linearly increasing  artery volume  up 100\% increase 
of $| \Omega^f_0|$ 
 since the inflow fluid volume over the simulation time $T$ is given by  $\int_0^T \int_{\Gamma_{in}}v_{0}\cdot\boldsymbol{n}_f \, dS\,  dt= |\Omega^f_0|$. 
Here $|\Omega^f_0|$ and $|\Gamma^f_{in}|$ are the volume of the sealed artery and the area of the inlet cross section, respectively and $\boldsymbol{n}_f$ is the outside normal vector.
Moreover, a smooth non-decreasing ramp function $0 \leq r(t)\leq 1$ is applied connecting the zero value with $v_0$  on a short  time interval to guarantee the consistency  of inlet boundary and initial condition and to develop the fluid movement from the initial state.

    \subsubsection*{Solid deformation subproblem}

The mechanics of the arterial wall is described by its deformation $\boldsymbol{d}:=\boldsymbol{x}(x,t)-X$ and can be modeled under some assumptions by the constitutive law for linear elastic, isotropic  material with following governing equations in the reference solid material domain $\Omega_0^s$,  
\begin{equation}\label{eq:soliddeformation}
    \begin{split}
        \rho_s \frac{\partial^2 \boldsymbol{d}}{\partial t^2} &= \nabla \cdot \boldsymbol{P}^T,\quad\text{ in } \Omega^s_0,\\
        \text{with}\quad
        \boldsymbol{P}&=J\left(\boldsymbol{C:\varepsilon}\right)\boldsymbol{F}^{-T},
    \end{split}
\end{equation}
where $\boldsymbol{F}=\frac{\partial x}{\partial X}, \ x \in\Omega^s_t, \  X \in \Omega^s_0 $, is the deformation gradient, 
$\boldsymbol{P}$ is the first Piola-Kirchhoff tensor and ${J=det\boldsymbol{F}}$.
The elastic strain tensor $\boldsymbol{\varepsilon}$ is given by the Green-Lagrange strain, ${\boldsymbol{\varepsilon}=\frac1{2}(\boldsymbol{F}^T\boldsymbol{F}-I)=\tfrac{1}{2}\left(\nabla\boldsymbol{d}+(\nabla\boldsymbol{d})^T+\nabla\boldsymbol{d}\left(\nabla\boldsymbol{d}\right)^T\right)}$. 
In the definition of $\boldsymbol{P}$ "$\boldsymbol{:}$" represents the tensor product.
The elasticity tensor 
\begin{equation}\label{def:elasticitytensor}
    \boldsymbol{C}=\boldsymbol{C}(E,\nu)
\end{equation}
depends on the material properties such as Young's modulus of elasticity ($E$) and the Poisson's ratio ($\nu$), the latter modeling the compressibility of the material.
In this contribution we  generalize the linear elasticity assumption by  varying Young's modulus dependent on strain by setting $E \equiv E(|\boldsymbol \varepsilon|)$. The dependency of $E$ on a strain metric $|\varepsilon|$ is observed in laboratory tensile tests of chosen material for fabricated CA-phantom, to  reproduce the experimentally observed inflation behavior of the phantom, 
the details on this modeling are given in section \ref{subsec:young_moduli}. Another source for the representation of Young's modulus is clinical data from \cite{faturechi2019mechanical}, which is investigated in section \ref{subsec:Invivo}.

The boundary of the solid domain representing the arterial wall is divided into three parts, $\Gamma^{s}=\Gamma^{s}_{in}\cup\Gamma^{s}_{ext}\cup\Gamma^0_{fsi}$.
For the inlet surface originating from the artificial clipping of the geometry,  $\Gamma^{s}_{in}$, the displacement is set to $\boldsymbol{d}=\boldsymbol{0}$. 
On the external surface $\Gamma^{s}_{ext}$, which is not in contact with the fluid,  the free stress boundary condition with no loads is imposed  $\boldsymbol P^T \boldsymbol{n}_s=0 $.   
%\todo{Tristan aus  COMSOL}
 \\

\subsubsection*{Fluid-structure interaction:}

The continuity of velocities and forces at the fluid-structure interface layer is maintained by the coupling conditions enforcing the balance of the velocities and stresses of the fluid and the solid material;
\begin{equation} \label{eq:couplings}
    \begin{split}
        \frac{\partial \boldsymbol{d}}{\partial t} = \boldsymbol{\tilde u}= \boldsymbol{w}\quad \text{and} \quad 
        J\tilde{\boldsymbol{T}}_f\,\boldsymbol{n}_f= -\boldsymbol{P}^T \boldsymbol{n}_s \quad \text{on }\Gamma^0_{fsi},
    \end{split}
\end{equation}
where $\boldsymbol{\tilde{u}}, \, \tilde{\boldsymbol{T}}_f$ are the velocity and Cauchy stress tensor of the fluid transformed to the reference fluid-solid interface, $\Gamma^0_{fsi}$.

\noindent
The coupled FSI problem (\ref{eq:ns}) -- (\ref{eq:couplings}) is numerically solved with the Finite-Element based software Comsol Multiphysics, Version 6.4 (COMSOL Inc., Stockholm, Sweden) utilizing the Structural Mechanics Module with monolithic (fully coupled) FSI approach. Linear P1-P1 elements were chosen for fluid and  quadratic serendipity elements were chosen for the discretization of solid deformation. 
The computational mesh is generated 
including several thin boundary layers at the vessel wall; the simulations for the compliance validation in section \ref{sec:sealedcarotis} are performed on a mesh consisting of almost 134 
tsd. tetrahedral elements for the fluid with  almost 41 tsd. prism-elements for two boundary layers
and using  ca 106 tsd. tetrahedral elements for the solid arterial wall discretization.
The numerical convergence of the solution with respect to the mesh refinement is validated in section \ref{sec:convergencestudy},  confirming the convergence of the obtained numerical solution  to a high-resolution reference solution.

\subsection{Experimental setup and stress-strain relation} 
 \label{subsec:laboratorysetup}

Numerical  modeling  of the carotid  artery volume change using above described FSI model is validated with laboratory measurements. Our previous validation results for sealed artery simulations compared with measurements on CA phantoms can be found in our work  \cite {ShiRiHu24}.
In this section, we summarize the attributes of the modeling setup that are relevant for investigating realistic stress-strain relations.

The physical artery is reconstructed using a silicon-like and a hydrogel material with similar mechanical properties to vessel tissue, which are also applied for fabrication of the  vessel implants PVA hydrogel and Dragon Skin material (DS) \cite{wang2021poly}. 
 The utility of such silicons for in vitro validation studies of humanoid arteries is further confirmed in   \cite{doyle2009experimental,souza2024experimental}.
The fabrication of artery phantoms is  based on so-called molding process. 
A deeper insight in the methods of geometric  reconstruction  of the physical artery phantom modeling is already given in 
 \cite {ShiRiHu24}. To determine the realistic  stress–strain relationship
% }
\begin{wrapfigure}{r}{0.5\linewidth}
    \centering
    \vspace{-0.3cm}
    \includegraphics[trim=0 460 40 50, clip,width=0.48\textwidth]{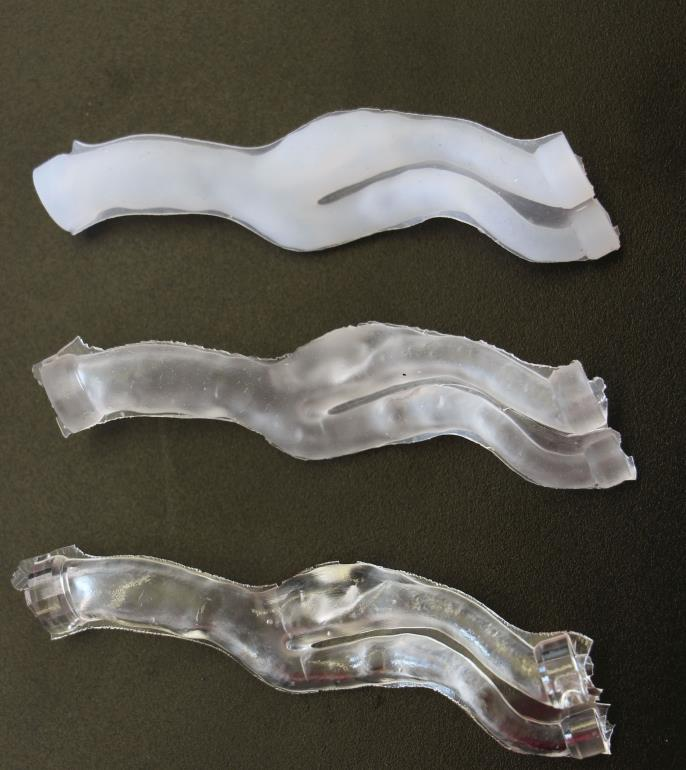}
    \caption{ Artery phantom from Dragon Skin}
    \label{fig:silicon}
    \vspace{-0.4cm}
\end{wrapfigure}
  for the elastic behavior  of the silicon artery  phantoms, experimental tensile  testing  to reconstruct the stress-strain relation was performed for cylindrical specimens of the material. 
These results are used to parametrize the  Young's elasticity module $E(|{\varepsilon}|)$  in the  elasticity tensor (\ref{def:elasticitytensor}).
To validate the measurement error of tensile tests,  five carotid  artery phantoms  were fabricated  by molding process.
In each of these artery reconstruction processes, additionally one material cylinder was fabricated with an identical material mixing ratio and supplied for tensile testing.
Measurements  were conducted three times with each specimen, totaling 15 measurements for each molded DS material.
Fig. \ref{fig:tensile_ds}  presents the measured stress curves $\sigma({\varepsilon})$ for DS probes as a function of uniaxial strain ${\varepsilon}$ and confirms the nonlinear stress-strain relationship when larger strain intervals are considered  (${\varepsilon} > 20 \%$).
\begin{figure}[htbp]
    \captionsetup{justification=centering, margin=0.5cm}
    \centering
     \includegraphics[trim=0 0 0 0, clip,width=0.87\linewidth]{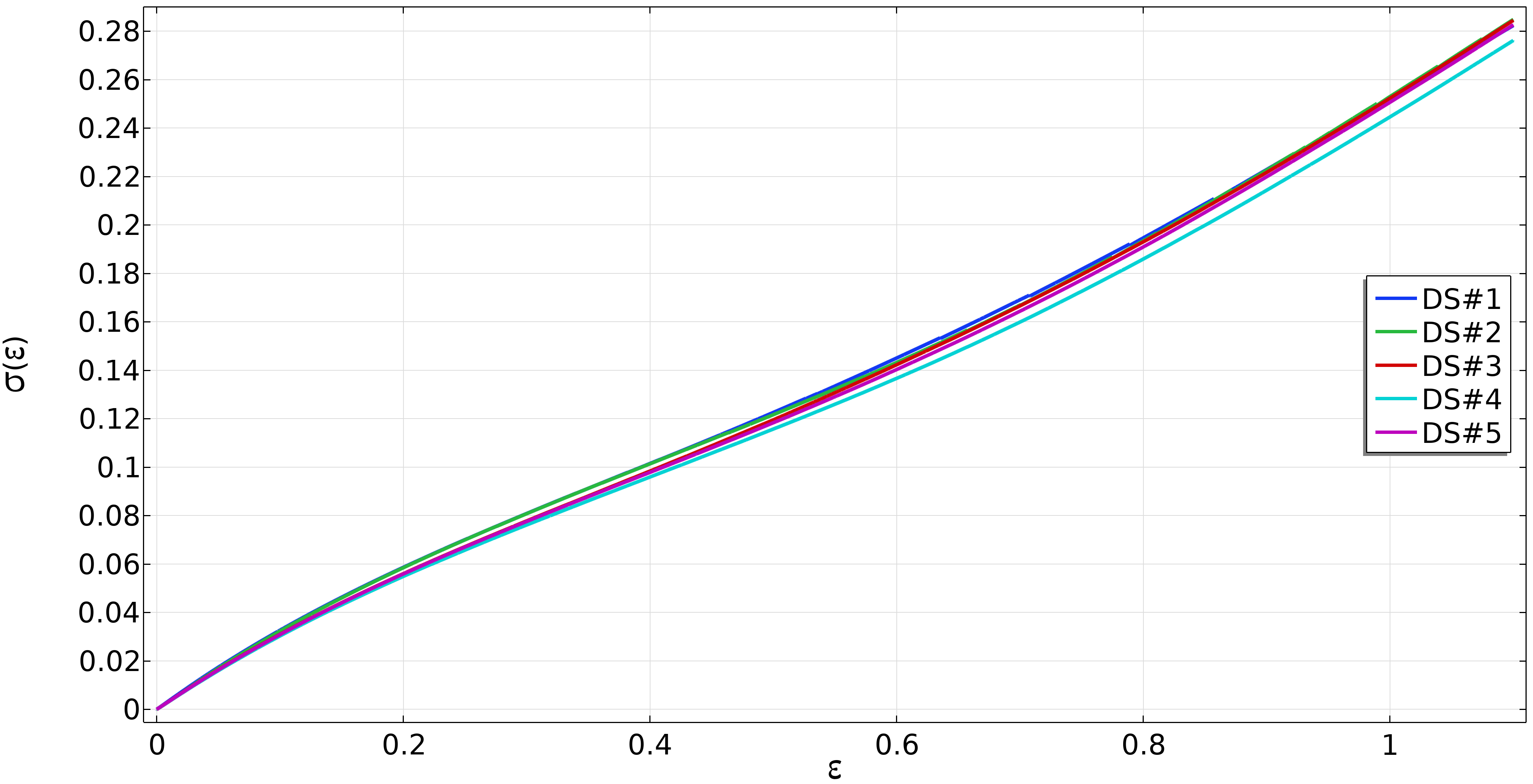}
     \caption{ Measured stress-strain behavior  of the uniaxial stress    $\sigma(\varepsilon)$ for five dragon skin (DS) probes, $\sigma$ in (MPa). Adapted from \cite[Fig. 8]{ShiRiHu24}\\.
     } 
     \label{fig:tensile_ds}
\end{figure}
\begin{figure}[htbp]
\captionsetup{justification=centering, margin=0.5cm}
    \centering
    \includegraphics[trim=0 0 0 0, clip,width=0.87\linewidth]{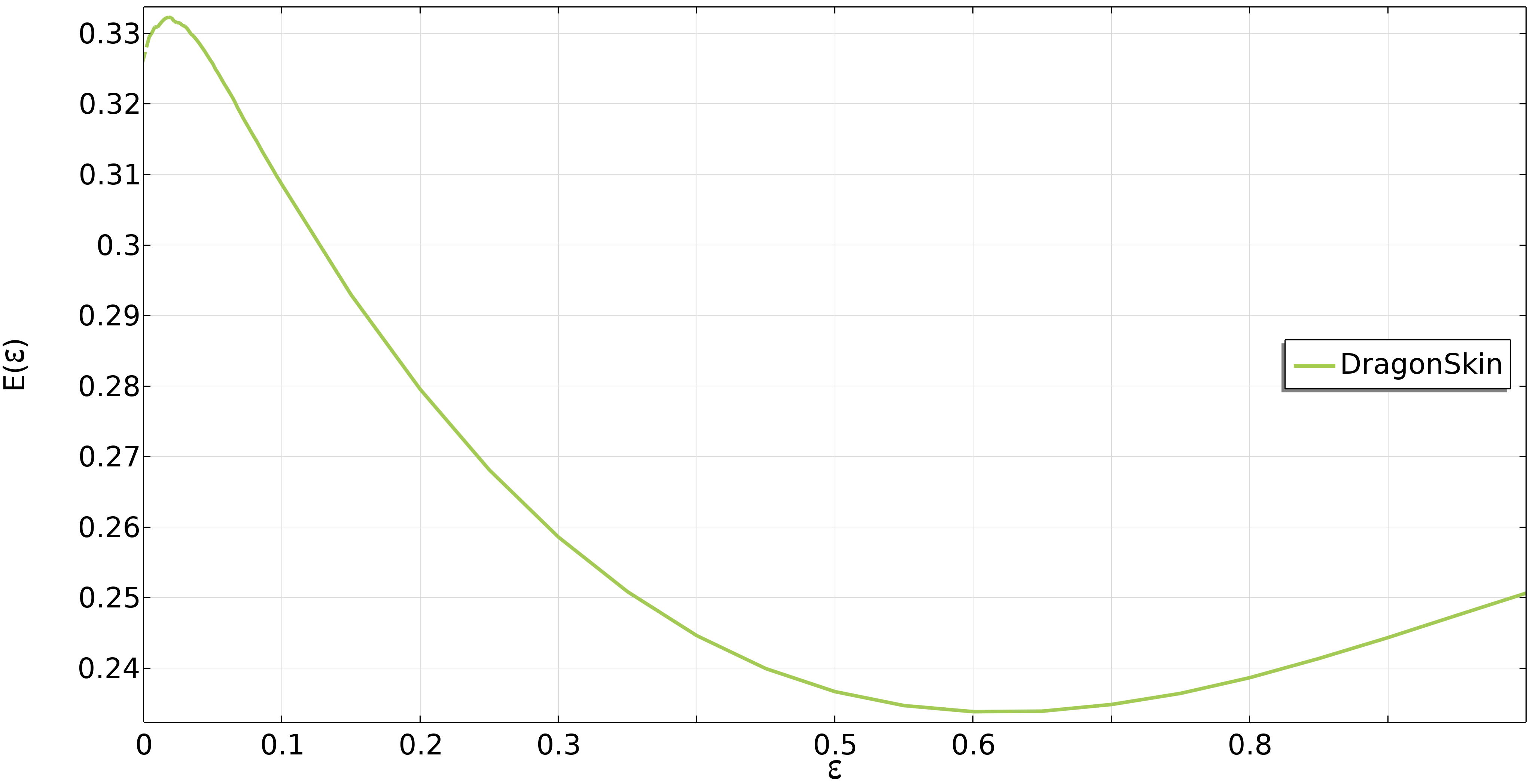}
     \caption{ Young's modulus function $E(\varepsilon)$ extracted as tangential slopes of the measured stress curves $\sigma(\varepsilon)$
     in (MPa). Adapted from \cite[Fig. 8]{ShiRiHu24}. }
     \label{fig:YoungsModulus}
 \end{figure}
The elasticity modulus  $E(\varepsilon)$ represents the 
(averaged) tangential slope of the measured stress curves $\sigma({\varepsilon})$  over the strain interval $(0,\varepsilon)$. 
It is obvious, that the initial slope  of the curves in Fig. \ref{fig:tensile_ds} overestimates the slope for larger strain ranges. 
The Young's modulus function $E(\varepsilon)$ 
is therefore constructed from stress curves averaged over material versions DS1 - DS5 in piece-vise manner  by extracting linear regression slopes $E(\varepsilon_k)$ of the stress curve $\sigma$ 
over strain intervals  $(0, \varepsilon_k), \ k=1,2 \ldots$. Strain intervals with stepwise increasing size,  $\varepsilon_{k+1}=\varepsilon_{k}+0.2\%$, 
up to $\varepsilon_{k+1}= 10\%$ and strain intervals 
with a longer increase of  5 $\%$ i.e., $ \varepsilon_{k+1}=\varepsilon_{k}+5\%$  for $\varepsilon_{k+1} >10\%$ are chosen in order to satisfactory represent the variations in tangential slope. 
The piece-vise linear fit  to the obtained  discrete values of $E(\varepsilon_k)$ presented in Fig. \ref{fig:YoungsModulus}, is considered to define the Young's modulus function $E(\varepsilon)$ 
in our generalized elasticity model with $\sigma(\varepsilon)=E(\varepsilon) \cdot  \varepsilon$. 
%
%
%-----------------------------------------
 \subsection{Young's modulus generalization}
\label{subsec:young_moduli} 
Based on the above described uniaxial nonlinear  stress-strain behavior we adapt the constitutive law  by a generalization of linear elasticity assumption $\sigma(\varepsilon)=E \cdot \varepsilon$ considering
  strain-dependency of the Young's modulus  and imposing $\sigma(\varepsilon)=E(\varepsilon) \cdot \varepsilon$. 

 In three dimensional case we replace the uniaxial strain with a specific scalar strain metric  $|\boldsymbol{\varepsilon}|$ and  employ the same dependency of Young's modulus on this scalar size, which results  in    
  strain-dependent elasticity tensor  $\boldsymbol C$, cf.  (\ref{def:elasticitytensor}) and nonlinear dependency of Cauchy stress tensor on strain,
\begin{equation}\label{eq:sigma_em}
    \boldsymbol{\sigma}:=\boldsymbol{\sigma} (\boldsymbol{\varepsilon}):=\boldsymbol{C}\left(E(|\boldsymbol{\varepsilon}|),\nu\right)\boldsymbol{:\varepsilon}
\end{equation}
with $|\boldsymbol{\varepsilon}|$ defined below.
The tensor product '$\boldsymbol{:}$' 
between the 4th order elasticity tensor $\boldsymbol{C}$ and the 2nd order strain tensor $\boldsymbol{\varepsilon}$ can be represented  in the Voigt notation for isotropic materials 
as a matrix - vector product  $\boldsymbol{D} \cdot \boldsymbol{\varepsilon}$ of the %a 3x3 
elasticity matrix $\boldsymbol{D}$ and the
%a 6x1 
strain vector 
$\boldsymbol{\varepsilon}=(\varepsilon_{xx}, \varepsilon_{yy}, \varepsilon_{zz}, \varepsilon_{xy}, \varepsilon_{xz}, \varepsilon_{yz})^T$,  with 
\begin{equation}\label{def:elasticitymatrixD}
\begin{split}
    \setlength\arraycolsep{2pt}
    \boldsymbol{D}=\frac{E(|\boldsymbol{\varepsilon}|)}{(1+\nu)(1-2\nu)}
    \begin{pmatrix}
        1-\nu & \nu & \nu &   \vdots  \\
        \nu & 1-\nu & \nu & \boldsymbol{0}  \\
        \nu & \nu & 1-\nu & \vdots \\
         \ldots & \boldsymbol{0}  &  \ldots & diag(\tfrac{1-2\nu}{2} )  \\
    \end{pmatrix},
\end{split}
\end{equation}
here $diag(\frac{1-2\nu)}{2}$ is a 3x3 matrix block.
The result of $\boldsymbol{D} \cdot \boldsymbol{\varepsilon}$
is a stress vector, which can be reordered into the symmetrical matrix $ \boldsymbol{\sigma}$ with entries $\sigma_{xx}, \sigma_{yy}, \sigma_{zz}, \sigma_{xy}, \sigma_{xz}, \sigma_{yz}$.

In 3D modeling, the strain metric which $E$ is dependent on is replaced by a specific scalar strain
$|\boldsymbol{\varepsilon}|$, 
cf. (\ref{eq:sigma_em}),  which accounts for both, the shear  as well as  the volumetric part of the  strain.  
Its exact definition can be identified based on the stress decomposition into the deviatoric and the volumetric part, $\boldsymbol{\sigma}=dev(\boldsymbol{\sigma})+\frac{trace(\boldsymbol{\sigma})}{3}$, resulting in the following formulation of the (linear) stress tensor
\begin{equation}\label{eq:stress_decomp}
   \boldsymbol{\sigma}=\boldsymbol{ C: \varepsilon }=\frac{E}{1+\nu} dev(\boldsymbol{\varepsilon}) +\frac{E }{ (1-2\nu)} \frac{\boldsymbol{\varepsilon}_{vol}}{3}\boldsymbol{I},
\end{equation}
with the volumetric strain $\boldsymbol{\varepsilon}_{vol}:=trace(\boldsymbol{\varepsilon})$ and  the deviatoric strain tensor $dev(\boldsymbol{\varepsilon}):=\boldsymbol{\varepsilon}-\frac{\boldsymbol{\varepsilon}_{vol}}{3}$.  
The volumetric strain is the first invariant of the strain tensor and  a scalar metric $\boldsymbol{\varepsilon}_{vol}=\varepsilon_{xx}+\varepsilon_{yy}+\varepsilon_{zz}$, which approximates  the volume change of a material cube by conserving its shape under the action of normal stresses $\sigma_{xx}, \sigma_{yy}, \sigma_{zz}$, 
\begin{equation}
    \frac{\Delta V}{V}=(1+\varepsilon_{xx})\cdot(1+\varepsilon_{yy})\cdot(1+\varepsilon_{zz})-1\approx\varepsilon_{xx}+\varepsilon_{yy}+\varepsilon_{zz},%=:\varepsilon_{vol}
\end{equation}
under the assumption of small strains.
The deviatoric strain 
is related to  the elastic shear strain $\gamma$ defined as 
\begin{equation}
    \gamma=\sqrt{2 \ dev(\boldsymbol{\varepsilon}): \ dev(\boldsymbol{\varepsilon})}=2\sqrt{J_2(dev ({\varepsilon}))},
\end{equation}
  where  $J_2$ is the second invariant of the strain deviator tensor
and describes the shape change 
%.e. the deformatioin 
of material cube under the action of shear (tangential) stresses $\sigma_{xy}, \sigma_{xz}, \sigma_{yz}$

\begin{figure}[htpb]
        \centering
        \includegraphics[trim={0cm 0cm 0cm 0cm},clip,width=0.75\linewidth, angle=0]{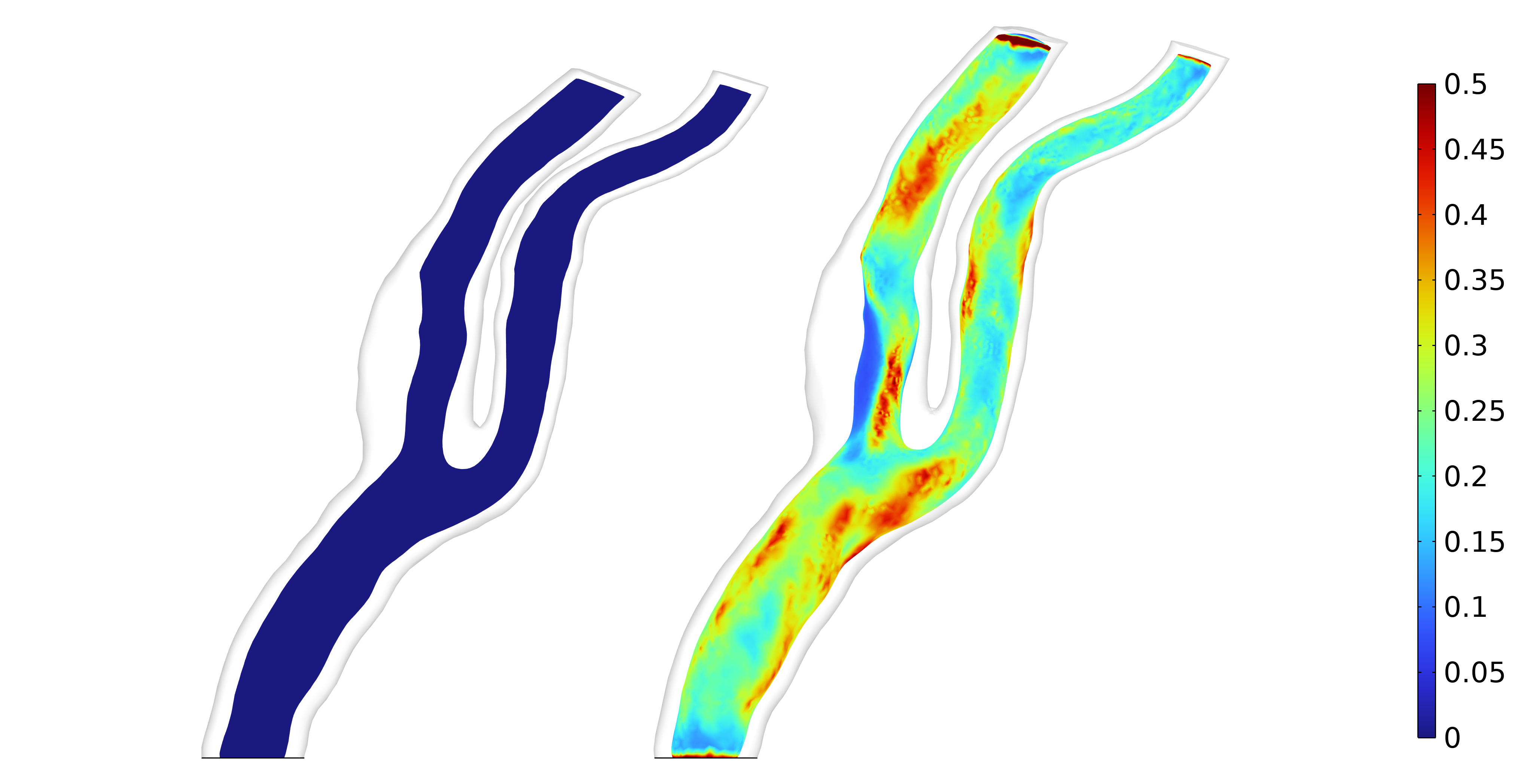}
       \hspace{0cm}
        \caption{ Surface distribution of the scalar strain  measure 
        $|\boldsymbol{\varepsilon}|$ (\ref{eq:em}), at initial state (left) and the time point of 45\% volume increase at  pressure $p= 148$ mmHg.}
        \label{fig:straindistribution}
\end{figure}

Based on the linear relation of the original Cauchy stress tensor (\ref{eq:stress_decomp})  with respect to  the strain tensor components $dev(\boldsymbol{\varepsilon}),  \ \boldsymbol{\varepsilon}_{vol}$  with proportionality constant  $E$, and based on the proposed generalization for the uniaxial stress 
%Young Moduli in uniaxial case 
$\sigma=E(\varepsilon) \cdot \varepsilon$,  we define the  strain metric 
$|\boldsymbol{\varepsilon}|$ 
%$\varepsilon_m$ 
accounting for both deviatoric as well as volumetric strain as follows
\begin{equation}\label{eq:em}
    %E=E(|\boldsymbol{\varepsilon}|), \quad
    |\boldsymbol{\varepsilon}|:=\frac{\gamma}{\sqrt{2}(1+\nu)}+\frac{\boldsymbol{\varepsilon}_{vol}}{3}\frac1{1-2\nu}. 
\end{equation}
Consequently, the strain dependent Young's modulus  $E=E(|\boldsymbol{\varepsilon}|)$  is considered in (\ref{eq:stress_decomp}).
The Young's modulus function $E(|\boldsymbol{\varepsilon}|)$ follows the slope curve extracted from  tensile tests presented in Fig. \ref{fig:YoungsModulus} and is applied in our numerical simulations of vessel inflation. 
The  surface distribution of the  strain measure 
$|\boldsymbol{\varepsilon}|$ extracted from simulations at 45 $\%$ volume increase
is presented in Fig. \ref{fig:straindistribution}
and shows up to 45-50 $\%$ strain increase  obtained locally during the carotid artery inflation at 
pressure  of 148 mmHg. This results in a non-negligible variation of the local Young's modulus along the strain range achieved, which impacts the evaluated compliance parameters.

%-----------------------------------------------------------------------------------------------------------------------------------------------------------------------------------------------------------------------------------------------------
\section{Validation of compliance of sealed carotid artery} \label{sec:sealedcarotis}

The compliance behavior of the sealed carotid artery is validated based on  the amount of deformation relative to the pressure increase during the artery inflation process. This  is evaluated in means of volume and diameter change as well as artery volume change with respect to the pressure exerted. 
In what follows, we summarize and present an extension of the results  published in \cite{ShiRiHu24}.
The presented simulation results have been obtained 
for the fluid viscosity $\mu=0.00345$ Pa.s,  the fluid density $\rho_f=1000$ kg.m$^{-3}$ and 
the Poisson ratio $\nu=0.4$. 
In the purely linear elastic model an averaged Young  modulus $E_{aver}$ set to $0.31$ MPa is used.
$E_{aver}$ is obtained as the slope of the linear fit to the stress–strain curves presented in Fig. \ref{fig:tensile_ds}  over the strain range up to $20 \%$, where the upper bound of  $20 \%$  was obtained as surface average of the simulated scalar strain measure
$|\boldsymbol{\varepsilon}|$ reached at pressure value of 120 mmHg representing the  physiological maximum.

\begin{figure}[htbp]
    \captionsetup{justification=centering, margin=1.5cm}
    \centering
    \begin{tikzpicture}
        
        \node[anchor=south west,inner sep=0] (img) {\includegraphics[width=0.8\linewidth]{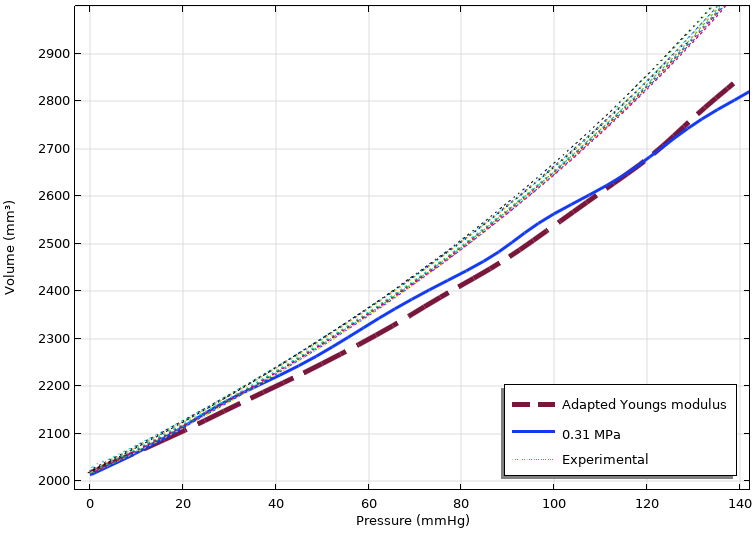}};

        \begin{scope}[x={(img.south east)}, y={(img.north west)}]
            
            % Hilfsgitter zum Finden der Position (nach dem Einstellen auskommentieren/löschen):
            % \draw[help lines, red, step=0.1] (0,0) grid (1,1); \foreach \x in {0,0.2,...,1} \node[red] at (\x,0) {\tiny \x}; \foreach \y in {0,0.2,...,1} \node[red] at (0,\y) {\tiny \y};
            
            % 3. Neuer Text mit weißem Hintergrund
            \node[fill=white, inner sep=2pt, font=\tiny, scale=0.8] at (0.85, 0.24) {Generalized Young's modulus};
            \node[fill=white, inner sep=2pt, font=\tiny, scale=0.8] at (0.78, 0.19) {0.31 MPa};
            \node[fill=white, inner sep=2pt, font=\tiny, scale=0.8] at (0.8, 0.14) {Experimental};
            
        \end{scope}
    \end{tikzpicture}
    \caption{ Pressure related artery volume increase, experimental (dotted) curves and simulated  (blue and red dashed) curves for linear and generalized elasticity model with Young's modulus function $E(|\boldsymbol{ \varepsilon}|)$ for Dragon Skin. Adapted from \cite[Fig 9.]{ShiRiHu24}}
    \label{fig:Silicon_Pressure_Volume}
\end{figure}
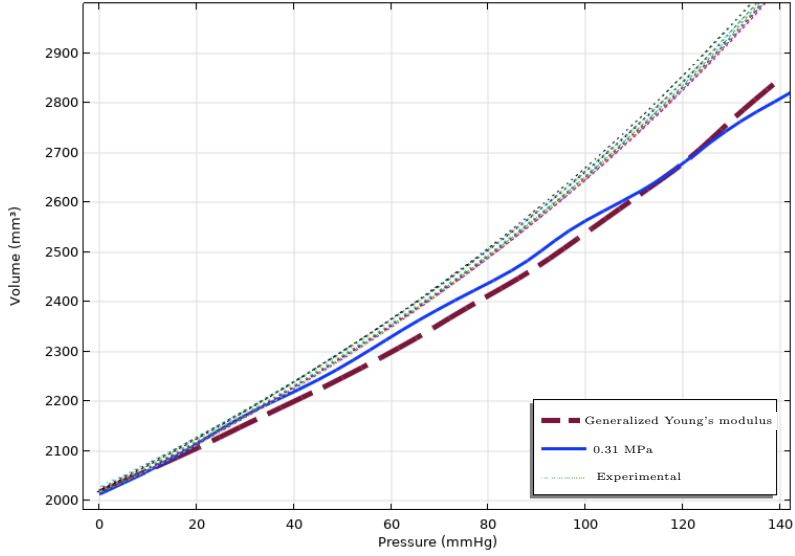

At first,  curves of volume change by an exerted pressure are qualitatively compared  in Fig. \ref{fig:Silicon_Pressure_Volume} for both measured and simulated values over a large pressure range to illustrate the volume increase. 
A clear improvement  with respect to the character of the stress-strain relationship has been obtained applying  the above approach of generalized Young's modulus dependent on strain  (\ref{eq:em}). An exponential relationship between fluid pressure and inflated artery volume for the generalized elasticity model can be observed in Fig. \ref{fig:Silicon_Pressure_Volume} (dark red dashed curve), whose exponential character is closer to the measured curves than the almost linear curve of the dependency of artery volume on pressure obtained for the pure linear elasticity model (blue curve). Here the impact of the local strain distribution observed on Fig. \ref{fig:straindistribution} on the Young's modulus justifies the differences in the pressure related volume increase.

Of particular biomedical interest  are the parameters of (relative) volumetric compliance $C$  and 
the (relative) diametric distensibility $D$ defined as
\begin{equation}\label{eq:compliance}
    C_{p_1-p_2}=\frac{V(p_2)-V(p_1)}{V(p_1)(p_2-p_1)}100 \%, \quad  D_{p_1-p_2}=\frac{D(p_2)-D(p_1)}{D(p_1)(p_2-p_1)}100 \%,
\end{equation}
with units [$\% / $mmHg]. 
They are obtained from the volume of inflated artery phantom, at some pressure  $V(p)$  or its diameter (at specific sites) $D(p)$, at two characteristic pressure values of the human cardiac cycle,  typical systolic/diastolic values  considered in literature are $p_1=74$  and $p_2=130$ mmHg or  $p_1=80$ and $p_2=120$ mmHg \cite{pascaner2015continuous}.
In this contribution we present the values of volumetric compliance and of the local diameter distensibility at later pressure range for Dragon Skin phantom see Table \ref{tab:1}, 
results for $70-130$ mmHg can be found in  \cite{ShiRiHu24}.  In our investigations the local diameter distensibility is evaluated at three specific diameter measurement sites, situated in the  common  carotid artery (CCA), in the external carotid artery (ECA)
and posterior to it (ICA),  see \cite[Fig. 5]{ShiRiHu24} for more details.
Photo-imagery is conducted  to measure the diameters in laboratory. 
For the simulated diameter values the artery cross-sections at the corresponding sites were projected into the measuring (y-z) plane, which is consistent with the laboratory setup.
\begin{table}[ht!]
\centering
\begin{tabular}{|l||c|c|c|c|}\hline 
 &$C_{80-120}$ & $D_{CCA}$ &  $D_{ECA}$ &$D_{ICA}$  \\\hline \hline
Measured values & 0.330\% & 0.168\% & 0.097\% & 0.151\%  \\\hline
Simulations, $E_{aver}=0.31$ MPa & 0.247\% & $-$ & $-$ &  $-$ \\\hline
%&$D_{ICAs}$& 0.069\%& $-$ & 0.129 \%
%Simulations, generalized $E(|\boldsymbol{\varepsilon}|)$   & 0.275\% & 0.144\% & 0.105\% & 0.155\% & 0.130\% \\\hline
Simulations, generalized $E(|\boldsymbol{\varepsilon}|)$   & 0.275\% & 0.144\% & 0.105\% & 0.160\%  \\\hline
\end{tabular}
\caption{Measured and simulated relative compliance ($C_{80-120}$) and distensibility ($D_{CCA/ICA/ECA}$) values obtained for pressures $80$ and $120$ mmHg.
The diameter distensibility is not evaluated for constant Young modulus. 
}
\label{tab:1}
\end{table}
Obviously, the obtained simulated relative volumetric compliance values are lower than the measured value, 
which can already be observed in the slopes of the secants to pressure-volume curves at considered pressures in Fig. \ref{fig:Silicon_Pressure_Volume}. However,  clear improvement  of the numerical  towards the measured compliance is obtained applying the generalized elasticity model with strain-dependent Young modulus $E(|\boldsymbol{\varepsilon}|)$, comparing to pure linear elastic model with  the constant value $E_{aver}$.
The diameter distensibility values are evaluated for the generalized elasticity model only:  at all sites
are they in a good agreement with the measured values.
{Obtained slight differences between simulation and measurements may explained by the fact that the uniaxial tensile test, on which the Young's modulus for simulation relies and the compliance experiment represent distinct mechanical loading scenarios. Vervenne et al. highlight the challenges arise for transition from tensile test to more complex geometrical 3D compliance scenario \cite{vervenne2025stretching}}. 
%

%-----------------------------------------------------------------------------------------------------------------------------------------------------------------------------------------------------------------------------------------------------
\section{Numerical convergence study for the sealed artery}
\label{sec:convergencestudy}

In what follows we validate the numerical convergence for  fluid velocity $\boldsymbol{d}$, pressure $p$ and deformation $\boldsymbol{u}$ at a fixed simulation time $t=0.11$ s, which corresponds to a 45\% volume increase of the sealed artery, i.e. the situation presented on Fig. \ref{fig:straindistribution} (right).
As a second, we present compliance quantities validates in section \ref{sec:sealedcarotis} at systolic-diastolic pressure range and compare them on consecutive refined meshes to demonstrate their 
convergence to a reference value obtained on highly resolved mesh.
The deformation problem has been solved using generalized elasticity model with generalized Young modulus dependent on  strain according to  (\ref{eq:sigma_em}), (\ref{eq:em}) using the Young's modulus function from Fig. \ref{fig:YoungsModulus} and setting $\nu=0.4$.

\begin{figure}[ht!]
    \captionsetup{justification=centering, margin=0.5cm}
    \begin{subfigure}{0.47\textwidth}
    \centering
    \includegraphics[width=\linewidth, height=6.3cm]{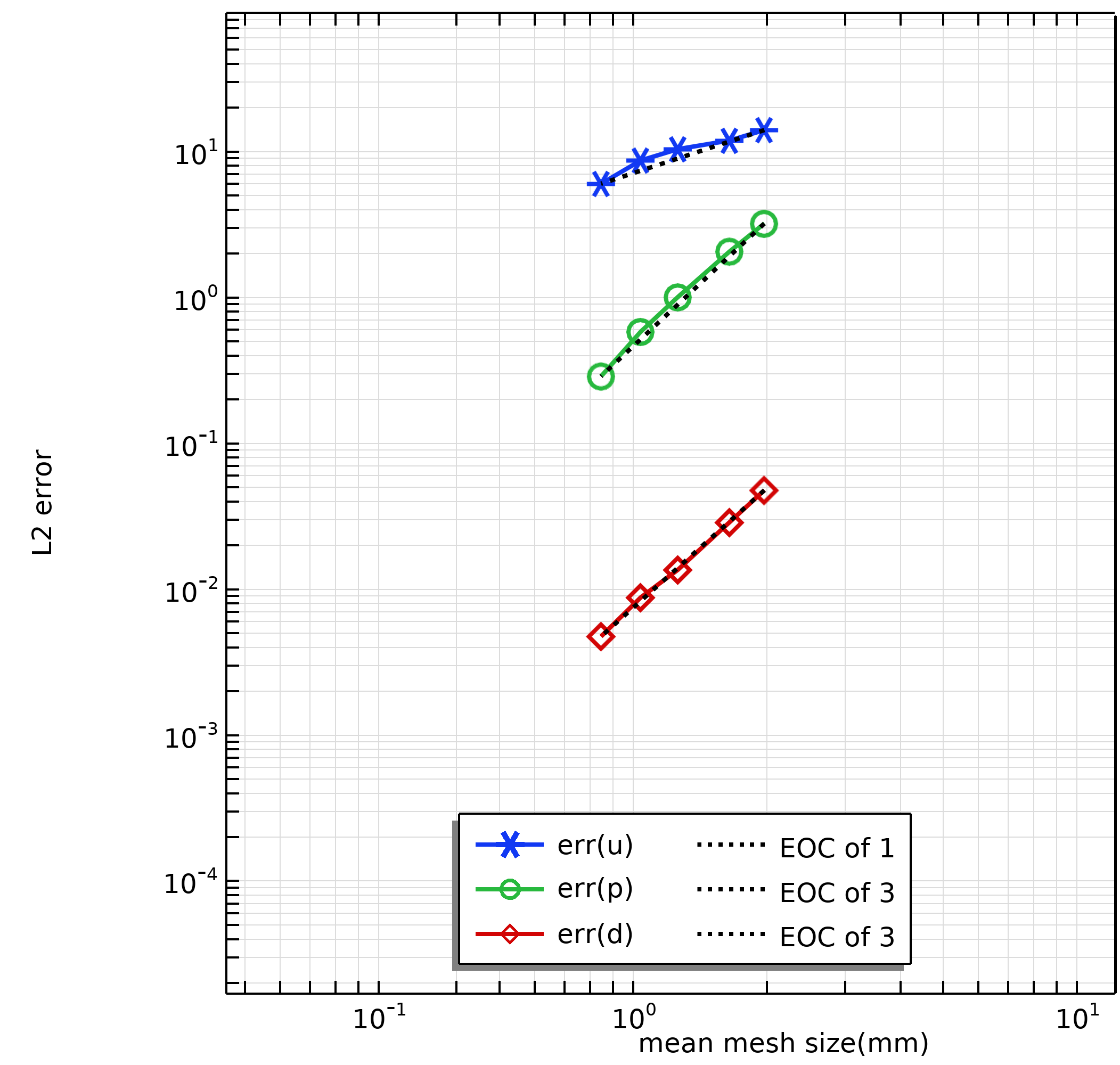}
    \caption{ Decrease of $L_2$ errors w.r. to mean mesh size  for velocity ($\boldsymbol{u}$), pressure ($p$) and wall deformation($\boldsymbol{d}$) at time $t=0.1s$.}
    \vspace{0.3cm}
    \label{fig:errorudp}
    \end{subfigure}
    \hspace{0.5cm}
    \begin{subfigure}{0.47\textwidth}
    \centering
    \includegraphics[ width=\linewidth, height=6.3cm]{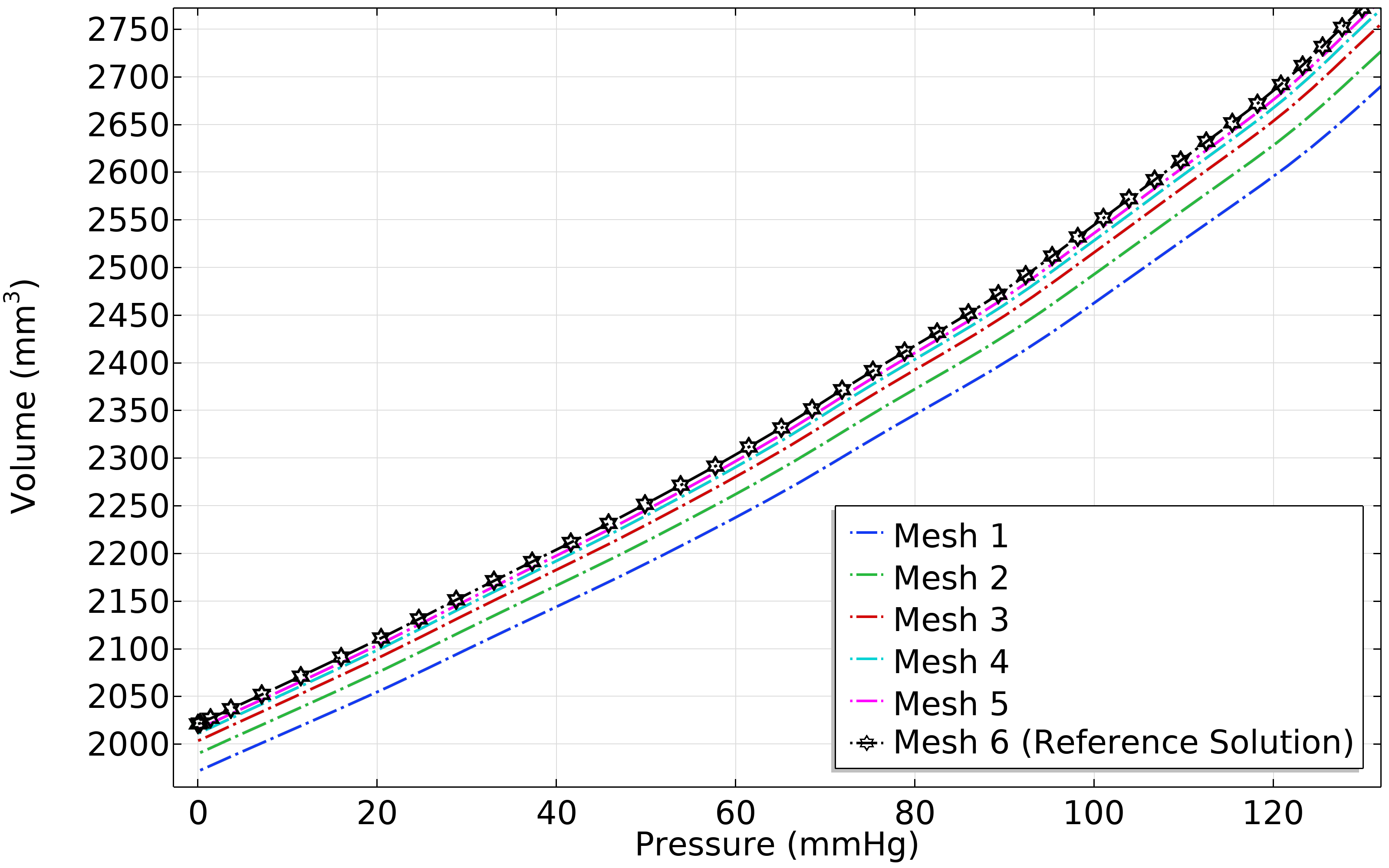}
    \caption{ Pressure related volume increase curves for the lowest up to the highest mesh resolution  show the convergence to the curve obtained on highly resolved mesh no. 6.}
    \label{fig:errorcurves}
    \vspace{-0.5cm}
    \end{subfigure}
    \vspace{0.2cm}
    \caption{ Convergence of $L_2$ errors for fluid flow and wall deformation (a) and of  pressure-induced volume increase curves (b).}
    \label{fig:meshconvergence}
\end{figure}

To demonstrate the convergence  with respect to the mesh refinement we use seven increasingly finer, unstructured meshes, 
the amount of elements roughly doubles in each consecutive mesh. The element number for each mesh is listed below in Table 
\ref{tab:physeval}.
The solution difference is evaluated by projecting the coarser solution onto the finest mesh, i.e., on  the reference mesh no. 6.
We use a weighted $L_2$-norm of the difference between the current and the reference mesh solution,
\begin{equation}
    \begin{split}\label{eq:errors}
    err(\boldsymbol{u}^{(i)})&:=\frac1{\sqrt{|\Omega^f|}}\|\boldsymbol{u}^{(i)}-\boldsymbol{u}^{(ref)}\|_{L^2(\Omega^f)}, \\
    err(p^{(i)})&:=\frac1{\sqrt{|\Omega^f|}}\|p^{(i)} - p^{(ref)}\|_{L^2(\Omega^f)},\\
    err(\boldsymbol{d}^{(i)})&:=\frac1{\sqrt{|\Omega^s|}}\|\boldsymbol{d}^{(i)}-\boldsymbol{d}^{(ref)}\|_{L^2(\Omega^s)}, \quad i=1, \ldots 5,\,
    \end{split}
\end{equation}
where $\boldsymbol{u}^{(i)}$, $p^{(i)}$ and $\boldsymbol{d}^{(i)}$ are the solutions at current $i$-th mesh. 
To identify the order of error decrease with respect to the mesh size,  the mean mesh size $\overline{h}$ is calculated as the mean value of the maximum inscribed spheres over all mesh element volumes.
The experimental order of convergence (EOC)
express the relative change of the errors of compared values (\ref{eq:errors}) with respect to the change of the mean mesh size in logarithmic scale,
\begin{equation} \label{eq:eoc}
    EOC (\boldsymbol{u}^{(i)})=\frac{\log_{10} (err(\boldsymbol{u}^{(i)}))-\log_{10}(err(\boldsymbol{u}^{(i+1)}))  }{\log_{10}{\bar{h}}^{(i)}-\log_{10} \, {\bar{h}}^{(i+1)}}
    =\log_{a^{(i)}}\left( \frac{err(\boldsymbol{u}^{(i)})}{err(\boldsymbol{u}^{(i+1)})}\right).
\end{equation}
Here $a^{(i)}:=\overline{h}^{(i)}/\overline{h}^{(i+1)}$ defines the decreasing factor of the mesh size,
ranging between 1.1961 and 1.3341 for meshes considered in this study.
%  1,2125, 1,2281,
The EOC (\ref{eq:eoc}) is analogously defined for $p$ and $\boldsymbol{d}$ in the analogy to the convergence order for the finite element method on a uniform mesh, where the numerical discretization error of the solution  $\boldsymbol{u}_h$ obtained with uniform mesh size $h$, $err(\boldsymbol{u}_h):=\|\boldsymbol{u}_h-\boldsymbol{u}_{ref}\|_{L_2}$ is expected to behave as $ C h^p, \ {C<\infty}$ with  a corresponding convergence order $p=log_a(err(\boldsymbol{u}^i)/err(\boldsymbol{u}^{i+1}))$ represented by the slope of the logarithmic error curves.
\begin{table}[ht]
\centering
\resizebox{\textwidth}{!}{%
\begin{tabular}{|c|c|c|c|c|c|c|c|}
\hline
\rule{0pt}{10pt}    
 Mesh & $\overline{h}$ (mm) & $err(\boldsymbol{u})$ (mm/s) & EOC($\boldsymbol{u})$ & $err(p)$ (mmHg) & EOC($p$) & $err(\boldsymbol{d}) $ (mm) & EOC($\boldsymbol{d}$) \\
\hline
\hline

1     & 1.9715 & 14.0489 & 0.9163 & 3.2146 & 2.4487  & 0.0477 & 2.7976  \\
2     & 1.6482 & 11.9031 & 0.5062 & 2.0642 & 2.6776 & 0.0288 & 2.7996 \\
3     & 1.2605 & 10.3922 & 0.9282 & 1.0067 & 2.8403 & 0.0136 & 2.2646 \\
4    & 1.0396 & 8.6904 & 1.7870 & 0.5824 & 3.4202 & 0.0088 & 2.9910 \\
5     & 0.8465 & 6.0196 &  & 0.2884 &  & 0.0048 &  
\end{tabular}}
\caption{Spatial discretization errors for velocity, pressure and deformation compared to reference mesh no. 6 and the corresponding EOCs at time $t=0.1s$}
\label{tab:numeval}
\end{table}
\begin{table}[!htbp]
\centering
\begin{tabular}{|c|c|c|c|c|c|c|c|}
\hline
\rule{0pt}{10pt}    
 Mesh & \# elements & $V_{120}$  & $err(V_{120})$ & $C_{80-120}$ & $err(C_{80-120})$\\
    &  & (mm$^3$) &  &  ($\%$/mmHg) &   \\
 \hline
 \hline
1 & 2.70$\times 10^4$ & 2595.45 &   3.35 \%     & 0.2663 &   3.75 \%          \\
2 & 4.71$\times 10^4$ & 2630.73 & 2.03 \% & 0.2723 & 1.58 \%  \\
3 & 1.05$\times 10^5$ & 2653.53 & 1.19 \% & 0.2728 & 1.39 \% \\
4 & 2.12$\times 10^5$ & 2667.32& 0.68 \% & 0.2744 & 0.83 \%  \\
5 & 4.41$\times 10^5$ & 2675.88  & 0.36 \% & 0.2755 & 0.42 \%  \\
% 6 & 1.26 $\times 10^6$ & 2683.86 & 0.29\% & 0.2765 & 0.35 \%  \\
6 & 1.77$\times 10^6$ & 2685.43 &  & 0.2767 & \\
\hline
\end{tabular}
\caption{Convergence of the  vessel volume $V$ at 120 mmHg, of  the volumetric compliance $C$ calculated at pressures 80 and 120 mmHg,  
including their relative errors $err(V), err(C)
$ compared to reference mesh no. 6, respectively.}
\label{tab:physeval}
\end{table}
The results of the error decay are presented in Table \ref{tab:numeval}
and the error curves are visualized in Fig. \ref{fig:errorudp}, 
indicating sub-linear convergence for the $L_2$-error  for fluid velocity, 
almost third order error decrease for the pressure and a order between 2 and 3 for the wall deformation. 
The expected second order $L_2$-convergence by applying linear finite elements for the velocity  is violated due to the inconsistency of the (spatially) constant boundary conditions of the inlet surfaces adjacent on its edges to the vessel wall equipped with non-slip condition. 
The atypically high EOC and the small errors of $\boldsymbol{u}$ and ${p}$ on mesh no. 5 can be explained by the proximity of the $5$-th mesh solution to the reference solution on mesh no. 7.

The accuracy of the  
compliance parameters is studied in means of the convergence of the vessel volume $V_{120}$ of the inflated phantom at systolic pressure $p=120$ mmHg, the convergence of the volumentric compliance $C_{80-120}$ (\ref{eq:compliance})  and of the difference of the pressure-volume curves (Fig. \ref{fig:Silicon_Pressure_Volume}), obtained for different mesh sizes.
The relative changes %
 of these quantities  are computed for each meshes $i$, $ \ i=1, \ldots 5$  and are evaluated w.r. to  the high fidelity reference mesh no. 6:
$
err(V_{120}^{(i)}):=\tfrac{|V_{120}^{(i)}-V_{120}^{(ref)}|}{V_{120}^{(ref)}}, \  err(C_{80-120}^{(i)}):=\tfrac{|C_{80-120}^{(i)}-C_{80-120}^{(ref)}|}{C_{80-120}^{(ref)}}.
$
The results of the mesh convergence study is presented in  Table \ref{tab:physeval} 
showing the decay of relative errors by a factor more than  ten within the chosen mesh sequence  for investigated compliance quantities.  
The convergence of the volume curves related to  the pressure is  presented in Fig. \ref{fig:errorcurves}.

We close this section a remark on expected numerical error in our volumetric compliance study presented in section \ref{sec:sealedcarotis}. The results for  the pressure related volume increase and of the volumetric compliance are computed on a mesh with a total about 280 tsd. elements, lying in the range between meshes 4 and 5 of the convergence study presented. 
Thus, from Table \ref{tab:physeval} one can deduce about  $0.83 \%$ to $0.41 \%$ expected relative error for the volumetric compliance  and about $0.68 \%$ to $0.36 \%$ relative  volume error
with respect to the solution on the highly resolved reference mesh.

%-----------------------------------------------------------------------------------------------------------------------------------------------------------------------------------------------------------------------------------------------------

\section{Compliance modeling for the open carotid artery }
\label{sec:openartery}

For the open compliant artery simulation,  the fluid-structure interaction model problem (\ref{eq:ns}),  (\ref{eq:soliddeformation}), with the coupling conditions (\ref{eq:couplings})  is extended by outlet boundary conditions at the artificial outlet cross sections $\Gamma_{out}$, see Fig. \ref{fig:realistic_geometry}. 
In absence of reliable clinical flow data for the daughter
CA-branches, the artificial outlets need to be accompanied by adequate boundary conditions that does not effect the flow patterns inside the artery volume. Windkessel models are a common class of physiologically reliable outlet boundary conditions \cite{taib2019wk4} which allow the modeling of downstream pulse pressure.  The serial four-element Windkessel model reads as follows.
\begin{equation}\label{eq:wk4}
     C L\frac{d^{2}Q}{dt^{2}}
    +\Big(CR_{p}+\frac{L}{R_{d}}\Big)\frac{dQ}{dt}
    +\Big(1+\frac{R_{p}}{R_{d}}\Big)Q
    =
    C\frac{dP}{dt}
    +\frac{P}{R_d}.
\end{equation}
Here, $R_p$ is the proximal, $R_d$ the distal resistance, C the downstream compliance and $L$ the inertance of the blood in the respective artery. $Q$ denotes  the volumetric flow rate and $P$ the outlets surface average pressure.
To model the compliance behavior in a physiological setting and to evaluate our generalized Youngs's modulus approach described in section \ref{subsec:young_moduli}, we simplify (\ref{eq:wk4})
by setting $C=0$ resulting in the following resistance-inertance condition (\ref{eq:resistance_inertia}),
\begin{equation}\label{eq:resistance_inertia} 
     P= L\frac{dQ}{dt}+R_{tot} Q.
\end{equation}
Here, $R_{tot}=(R_d+R_p)$ equals the total resistance of the downstream region.
Conceptually,   the total resistance $R_{tot}$ is defined as a constant relationship between the mean pressure $P_{mean}$ and the (mean) volumetric flow rate $Q_{mean}$ achieved
in the artery lumen $\Omega^f_t$, thus  $P_{mean}=Q_{mean}\cdot R_{tot}$
is directly proportional to the pressure $P_{mean}$ that needs to be maintained to achieve a certain volumetric flow rate $Q_{mean}$ \cite{pape2019physiologie}.
The total resistance of the branched artery $R_{tot}=P_{mean}/{Q_{mean}}$ , is further distributed into the the sub-branches of the arterial tree, which can be considered as parallel branches of  a part of the circulatory system.
Therefore, the total resistance is equal to the sum of the  
reciprocal values of the sub-branch resistances,
\begin{equation}\label{eq:resistances}
    \frac1{R_{tot}}=\frac1{R_{int}}+\sum_{i=1}^{N_{ext}}\frac1{R_{ext_i}},
\end{equation}
where  $R_{int},R_{ext_i}\,,i=1\dots N_{ext}$ denote the resistance values of the internal and  ($N_{ext}$) external branches of CA, respectively.
Having
specified the total resistance $R_{tot}$, 
the  key aspect 
is the distribution of flow volume between the internal and the external CA branches according to (\ref{eq:resistances}). 
According to  \cite{elwertowski2018importance}, the cranial supply role of the internal CA is maintained through various degrees of stenosis in the sinus bulb. 
We employ  this fact by  assuming a fixed   
 outflow volume ratio $\alpha$  between the internal CA and the external CA branches,  s.t., $$Q_{int}=\alpha Q_{mean}, \quad Q_{ext}=(1-\alpha)Q_{mean}, \ 0<\alpha <1.$$
Generally, the fixed division ratio holds for the volume flow averaged over a longer period of time, but represents a strong simplification with respect to the instantaneous dynamics.
Consequently, using the definition of particular
resistances of internal/external CA branches, $R_{int}={P_{mean}/}{Q_{int}}, \ R_{ext}={P_{mean}}/{Q_{ext}}$, applying the above flow division  ratio $\alpha$ for $Q_{int/ext}$, for the resistance value of the internal CA  and for the collective resistance of all external CA branches, one gets 

\begin{equation}
\label{eq:RinRex}
    R_{int}=\frac1{\alpha} R_{tot}\; \text{ and }\; R_{ext}=\frac1{(1-\alpha)}R_{tot}.
\end{equation}
Finally, we apply Murray's law \cite{FevRozzaMurray2021} to distribute the flow through the sub-branches of the ECA proportionally to their diameters, leading to the corresponding resistances being inversely proportional to the proportion of the cross-sectional outlet areas $|\Gamma_{ext_k}|$ of ECA sub-branches, 
\begin{equation}\label{eq:Rex}
    R_{ext_k}=\frac{\sum_{i=1}^{N_{ext}} |\Gamma_{ext_i}|}{|\Gamma_{ext_k}|}R_{ext}\text{, for }k=1\dots N_{ext}.
\end{equation}
\newline
The inertance $L_{i}$ for the i-th arterial segment on the other hand is defined  as a relationship between the pressure drop and the rate of change of volumetric flow $L_{i}={\Delta P}/\frac{dQ}{dt}$, describing how much pressure must be applied to accelerate or slow down the volumetric flow. $L_i$ can be approximated as mentioned in  \cite{stergiopulos1999inertanca} as
\begin{equation*}\label{eq:Li}
    L_{i} = \int_{0}^{l_{i}} \frac{\rho}{r_{i}^{^2}(x)\text{ }\pi}dx,
\end{equation*}
where $l_{i}$ is segment length, $\rho$ blood density and $r_{i}(x)$ the local radius of the arterial segment.
With respect to available data it is convincing to use a constant mean radius $\bar{r}_{i}$ and consequently the respective arterial segment  inertance becomes
\begin{equation}\label{eq:Li_mean}
    L_{i} = \frac{l_{i} \text{ }\rho}{\bar{r}_{i}^{^2}\text{ }\pi}.
\end{equation}
The total inertance for parallel arteries is obtained analog as for resistance: ${L_{tot}^{-1}}={L_{int}^{-1}}+L_{ext}^{-1}$. Here, if the total inertance is known, for instance from in vivo data, the segmental inertances are calculated by exploiting the ratio $k_{L}=\frac{L_{ext}}{L_{int}}=\frac{{l_{ext}}/{\bar{r}^2_{ext}}}{{l_{int}}/{\bar{r}^2_{int}}}$, s.t.,
\begin{align}\label{eq:Lseg}
    L_{ext}=L_{tot} (k_{L}+1), \ \ \ \ \ \
    L_{int}=L_{tot} \frac{1+k_{L}}{k_{L}}.
\end{align}
 In the absence of specific geometric data for the externa daugther branches (ECA's), we assume a uniform length $l_i$  for the distribution of $L_{ext}$ to the sub-branches.
By applying the aforementioned procedure for $L_{ext}$ and $L_{int}$ to $L_{ext_i}$,  under the assumption of uniform lengths $k_{L}$ becomes ${\bar{r}^2_{int}}/{\bar{r}^2_{ext}}$ and the distribution of inertance behave analogous to resistance in accordance with Murray's Law:

\begin{equation}\label{eq:Lex}
    L_{ext_k}=\frac{\sum_{i=1}^{N_{ext}} |\Gamma_{ext_i}|}{|\Gamma_{ext_k}|}L_{ext}\text{, for }k=1\dots N_{ext}.
\end{equation}
The above derived downstream resistances (\ref{eq:RinRex}), (\ref{eq:Rex}) and inertances (\ref{eq:Lseg}), (\ref{eq:Lex})   are implemented in (\ref{eq:resistance_inertia}) as the following implicitly formulated boundary condition at the outlet cross sections, i.e. prescribing
\begin{equation} \label{eq:pout}
    P_{i} = R_{i}\int_{\Gamma _i}\boldsymbol{u}\cdot \boldsymbol{n}_{f} dS + L_{i}\frac{d}{dt}\Big(\int_{\Gamma _i}\boldsymbol{u}\cdot \boldsymbol{n}_{f} dS\Big),
\end{equation}
where i $\in int,ext_{1},ext_{2},...,ext_{N_{ext}}$ denote the corresponding sub-branches of the arterial tree.

Finally, we consider a given time dependent inlet velocity $v_{in}(t)$ and zero deformation at inlet/outlet (clamped) surfaces of the open CA-artery and -tissue,
\begin{align}
    \boldsymbol{u}\cdot \boldsymbol{n}_f&={v}_{in}(t), \quad  \quad
    \boldsymbol{u}\cdot t_f=0 \quad  &\text{on }\ \Gamma^f_{in}\,,& \label{eq:init+boundCond1}\\
   % p&=-R_{s}\int_{\Gamma_{s}}\boldsymbol{u}(t)\cdot\boldsymbol{n}\,\text{dS}\,,\quad s\in\{int,ext_1,\dots,ext_3\} &\text{on }&\Gamma^f_{out},& \label{eq:init+boundCond2} \\
    \mathbf{d}&\equiv\mathbf{0} \quad &\text{on }\Gamma^s_{in}\cup\Gamma^s_{out},& \label{eq:init+boundCond3}
\end{align}
with $\boldsymbol{t}_f$ are the inlet surface tangential vectors.
%

%%%%%%%%%%%%%%%%%%%%%%%%%%%%%%%%%%%%%%%%%%%%%%%%%%%%%%%%%%%
\subsection{ Validation of the thin-walled CA model}
%%%%%%%%%%%%%%%%%%%%%%%%%%%%%%%%%%%%%%%%%%%%%%%%%%%%%%%%%%%

% Numerical setup
In the following  a realistic thin-walled artery model cropped with one inlet and four open outlet cross-sections,  as depicted in Fig. \ref{fig:geometries} is considered.
As the ratio of vessel wall tissue to flow lumen volume in this realistic geometry differs from the capped/sealed CA model used for compliance validation in section \ref{sec:sealedcarotis}, the Young's modulus is scaled to compensate for the reduction in wall volume following Laplace's law \cite{fung1993biomechanics}.
A multiplicative factor $k$ is introduced into the Young's modulus for the thin-walled tissue, defined as  $E_{th}(|\boldsymbol{\varepsilon}|)=k \cdot  E(|\boldsymbol{\varepsilon}|)$. This factor is calibrated to ensure that the pressure-induced volume change of the realistic artery tree with sealed outlets (Fig. \ref{fig:realistic_geometry}) matches that of the thick-walled sealed artery (Fig. \ref{fig:cropped_geometrie}), as validated in section \ref{sec:sealedcarotis}. In this manner a value of  $k=2.45$ is determined, which preserves the simulated pressure-volume behavior of the thick-walled reference for the thin-walled model under the same pressure loading, compare the dark red dashed curve in Fig. \ref{fig:Silicon_Pressure_Volume}.
Subsequently we consider the scaled Young modulus function $E_{th}(|\boldsymbol{\varepsilon}|)=2.45 \cdot  E(|\boldsymbol{\varepsilon}|)$  for the cropped patient-specific thin artery with four open outlets. 
At the inlet cross section we consider the time-dependent (spatially constant) normal inlet velocity  ${v}_{in}(t)\ =Q_{CA}(t)/|\Gamma_{in}^f|$ in (\ref{eq:init+boundCond1}),
%with a distinct flow velocity function $\boldsymbol{v}_{in}(t) $ 
that  produces typical volume flow rate  wave of the carotid artery $ Q_{CA}(t)$ presented in Fig. \ref{fig:inflow}, reflecting   systolic and diastolic peaks and a cardiac output of $432$ ml/min \cite{antiga2002patient,perktold1995computer} . Consequently the mean volumetric flow rate is validated,
\begin{equation}\label{eq:meanflowvolume}
    Q_{mean}=\frac1{T}\int_{0}^{T}Q_{CA}(t)\,dt=\frac1{T} \int_{0}^{T}
    {v}_{in}(t) \cdot \Gamma_{in}^f
    %\boldsymbol{v}_{phys}(t)\mathbf{n}
    \,\,dt\,=\text{7,2 ml/s}.
\end{equation}
The mean pressure together with $Q_{mean}$ drives the flow and defines the total resistance in carotid artery $R_{tot}=Q_{mean}/P_{mean}$, where $P_{mean}$ can be approximated by  \cite{levick2013introduction}, 
\begin{equation}\label{eq:meanpres}
    P_{mean}=\frac{P_{sys}+2\cdot P_{dias}}{3},
\end{equation} 
with $P_{sys},\ P_{dias}$  the reference systolic and diastolic pressure values. 
We consider $P_{sys}=120 \text{ mmHg and }P_{dias}=80$ mmHg in this study
resulting  in $P_{mean}=93,33$ mmHg  and  the total resistance of the carotid artery  $R_{tot}\approx 12.96 $ based on (\ref{eq:meanflowvolume}), (\ref{eq:meanpres}).
The resistances $R_{int}, \ R_{ext}$ of the interna and externa CA branches considered in (\ref{eq:RinRex})  are then determined based on constant flow division factor set to $\alpha=70\%$ dictating constant $70\%:30\%$ flow division ratio between the internal and external CA, respectively. 
 This is roughly in agreement with typical,  averaged over the cardiac cycle MRI-based flow division ratios for the CA,  e.g.,  $75 \%:25\%$ obtained in \cite{MILNER1998} or $64\%:35\%$ \cite{outflow_murray2010} for low stenosis degree.
 \begin{table}
    \centering
    \begin{tabular}{|c|c|c|}
    \hline
         & $R_i$ [$\frac{mmHg \cdot s}{ml}$] & $L_i$[$\frac{mmHg \cdot s^2}{ml}$]   \\
         \hline
         Interna& 18.502 &0.202215   \\
         \hline
         Externa&43.172 &0.09584    \\
         \hline
         Total&12.96 & 0.064         \\
         \hline
    \end{tabular}
    \caption{Resistance and inertance parameter values for ICA and  ECA}
    \label{tab:paramters}
\end{table}
 In  \cite{Marshall_2004}, %
 the time-averaged ICA/CCA and ECA/CCA ratios for a cohort of patients were obtained as   $70\%$ and $26\%$, which correspond approximately with our approach for fixed division ratio. See also our later comments on results and comparisons presented in Fig. \ref{fig:flow_splitting} to this point.
The total inertance  $L_{tot}=0.064 \frac{mmHg \text{ }s^2}{ml}$ for CCA is considered according to \cite{paisal2019analysis}. The inertances $L_{int}, L_{ext}$ for each branch are calculated as described by (\ref{eq:Lseg}) using the inertance ratio factor $k_{L} \approx 0.5$ based on mean radii $\bar{r}_{int}$,$\bar{r}_{ext}$ and lengths $l_{int}$, $l_{ext}$  for ICA and ECA found e.g., in \cite{baz2021morphometry,charles2021origin,choudhry2016vascular,sasikumar2023morphometric}. 
For the subbranches of ECA we apply  Murray's law, c.f.,  (\ref{eq:Rex}),(\ref{eq:Lex}) to define  $R_{ext_i}, L_{ext_i}, i=1,2,3$.
Finally, the obtained resistance and inertance values,  as presented in Table \ref{tab:paramters} are considered in boundary condition (\ref{eq:pout}) and implemented in explicit way.
\begin{figure}[ht!]
    \captionsetup{justification=centering, margin=1.5cm}
     \centering
     \begin{tikzpicture}
        \tikzset{
            myarrow/.style={-{Stealth[scale=1.2]}, thick, draw=black, rounded corners=3pt},
            labelstyle/.style={ fill=white, fill opacity=0.85, text opacity=1, inner sep=2pt}
        }
    % 1. Hauptbild platzieren
    \node[anchor=south west, inner sep=0] (main) {\includegraphics[width=0.75\linewidth]{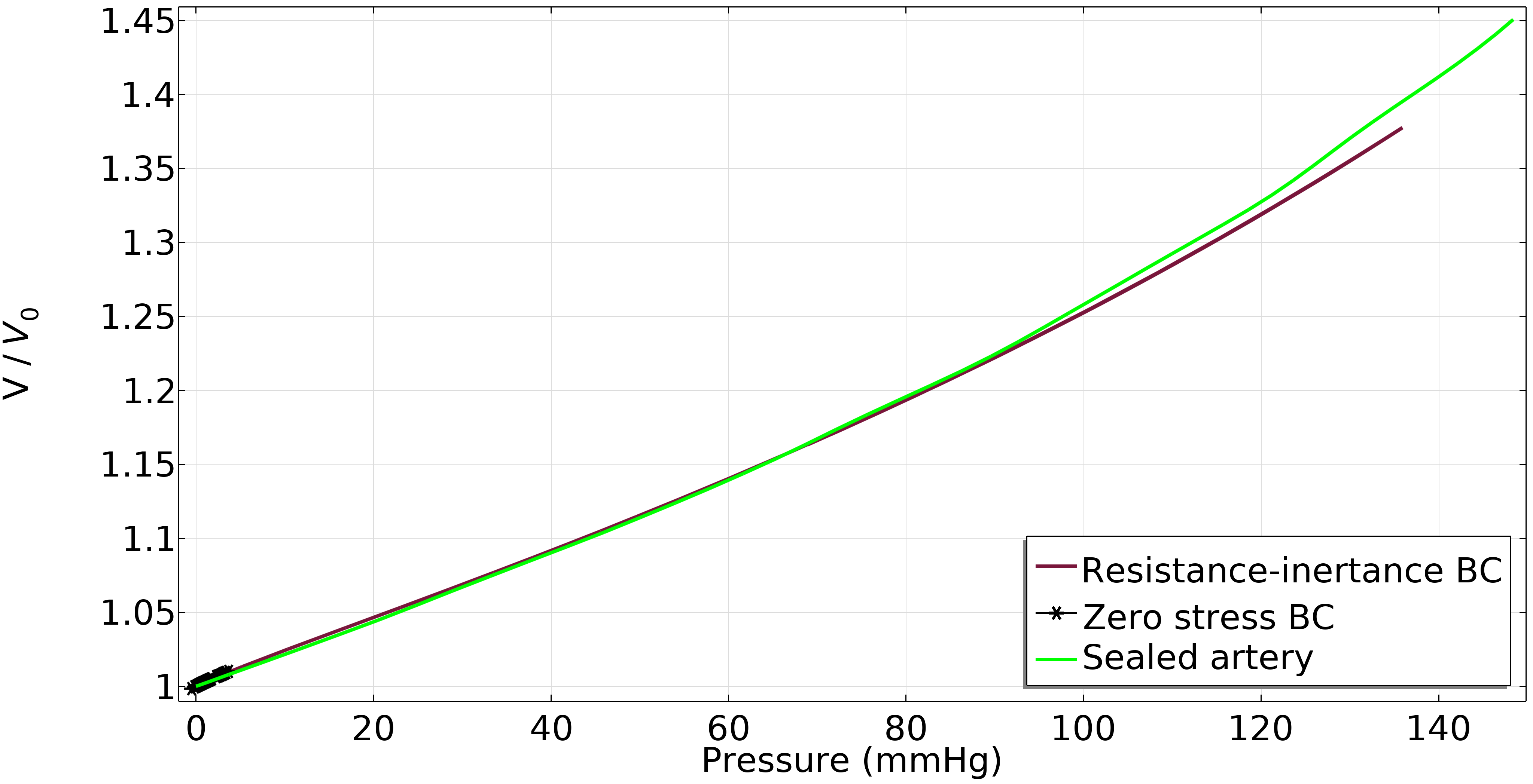}};

    \begin{scope}[x={(main.south east)}, y={(main.north west)}]
        
        \node[inner sep=0, draw=white, line width=1pt] at (0.3, 0.7) {
            \includegraphics[trim= 4cm 1.7cm 7cm  11cm, clip,width=3cm]{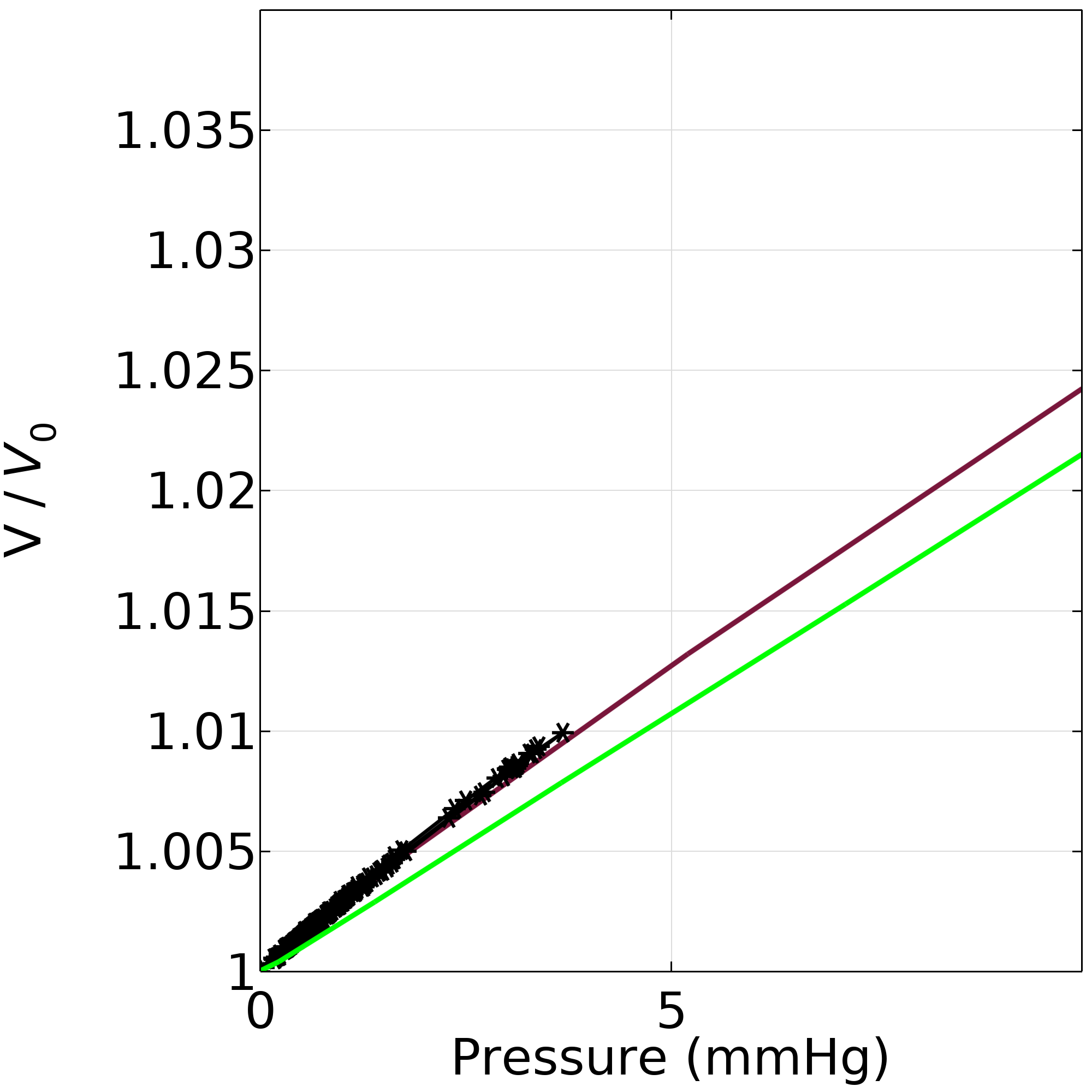}
        };
         \draw[black, thick, fill=gray, fill opacity=0] 
                    (0.14, 0.5) rectangle (0.46, 0.92);

        \draw[black, thick, fill=gray, fill opacity=0] 
                    (0.12, 0.11) rectangle (0.16, 0.16);
                    
        \draw[black, thick] (0.12,0.16) -- (0.14,0.5);
        \draw[black, thick] (0.16,0.16) -- (0.46,0.5);

    \end{scope}
  \end{tikzpicture}
     \caption{  Validation of the pressure-induced relative volume change for the open artery model using resistance-inertance boundary conditions. Results are compared to the sealed artery reference and to zero-stress BC.}
     \label{fig:pv_modelvalidation}
 \end{figure}
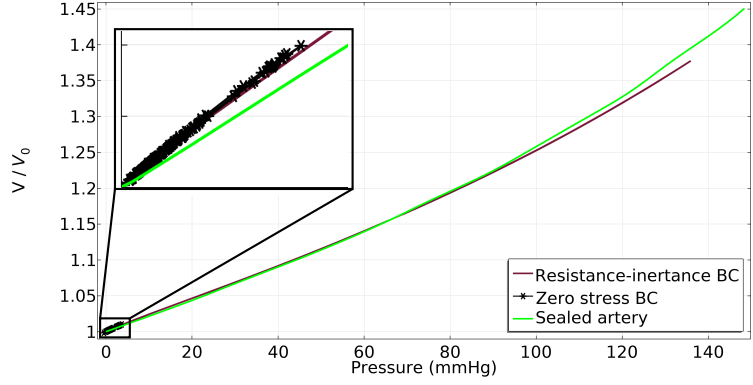

The model setup is validated with respect to the pressure-volume behavior, expressed as relative volume $\frac{V}{V_0}$, shown in Fig. \ref{fig:pv_modelvalidation} for generalized Young's modulus $E(|\boldsymbol{\varepsilon}|)$ cf. (\ref{eq:sigma_em}),(\ref{eq:em}), where 
$V_0$
denotes the volume at zero (starting) pressure. To this end, the pressure-volume curves of the open artery case employing the resistance-inertance boundary condition (\ref{eq:pout}) are compared against two reference scenarios: the in vitro validated sealed artery simulation with strain dependency of Young's modulus as presented in section \ref{sec:sealedcarotis}, and the open artery case implementing the commonly applied zero-stress boundary condition ($\boldsymbol{T}_f\boldsymbol{n}_f=0$) for free outflow.
\begin{figure}[ht!]
    \captionsetup{justification=centering, margin=0.5cm}
    \centering
     \includegraphics[width=0.49\linewidth,height=5cm]{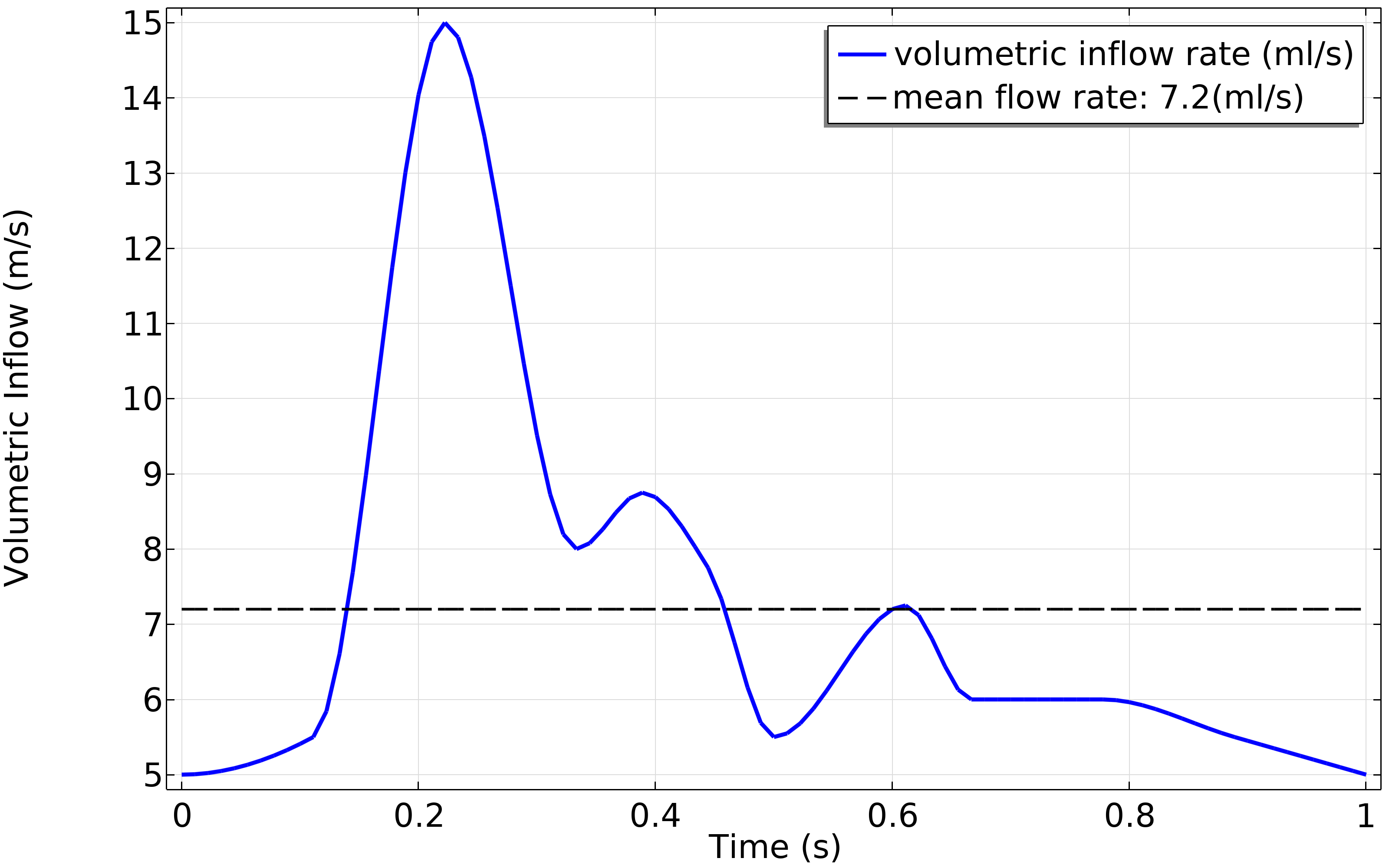}
    \hfill
    \includegraphics[width=0.49\linewidth,height=5cm]{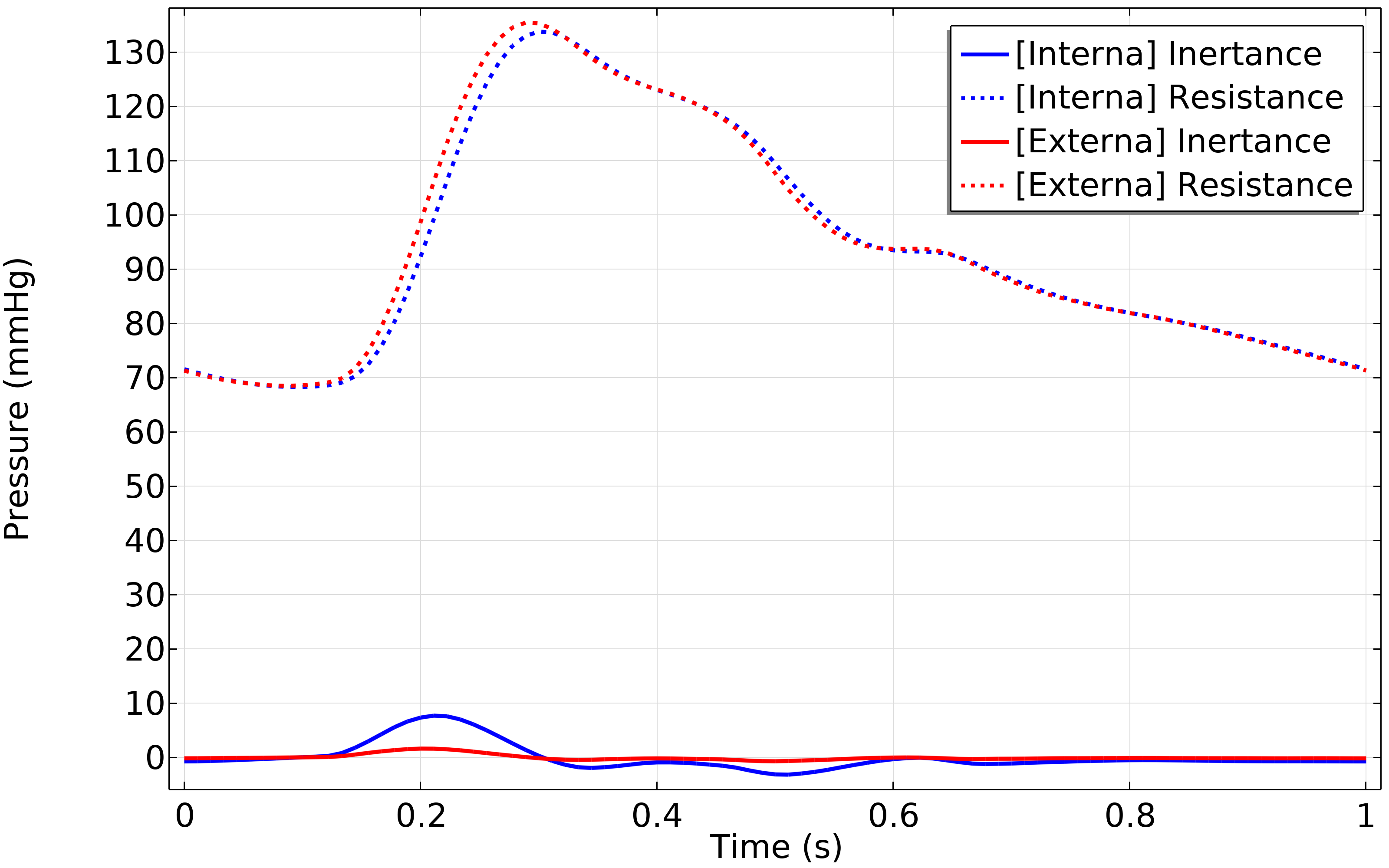}\\[1ex]
    \caption{ Typical volumetric flow rate  $Q_{CA}(t)$ \cite{antiga2002patient,perktold1995computer,ProRiHu23} for inflow BC at $\Gamma^f_{in}$ (left), and temporal evolution of the inertance and resistance components of the pressure condition (\ref{eq:pout}) for the ICA and ECA over one normalized cardiac cycle (right).}
    \label{fig:inflow}
\end{figure}
An almost perfect match is observed between the resistance-inertance boundary condition (in red) and the previous sealed artery simulation (in green). This agreement confirms the successful transfer of the volume behavior between the sealed (in vitro validated) and open configurations for corresponding pressure action. In contrast, the FSI simulation utilizing the zero-stress boundary condition (black cross) fails to generate a physiologically meaningful pressure buildup and, consequently, does not exhibit an significant inflation. This shows that the zero-stress approach is unsuitable for the simulation of arterial fluid and  wall mechanics and confirming clearly a strong impact of the outlet boundary condition on the compliance behavior, which is properly reconstructed by application of the proposed resistance-inertance BC.
    \begin{figure}[ht!]
    \captionsetup{justification=centering, margin=0.5cm}
            \centering
        \includegraphics[width=0.75\linewidth,height=5cm]{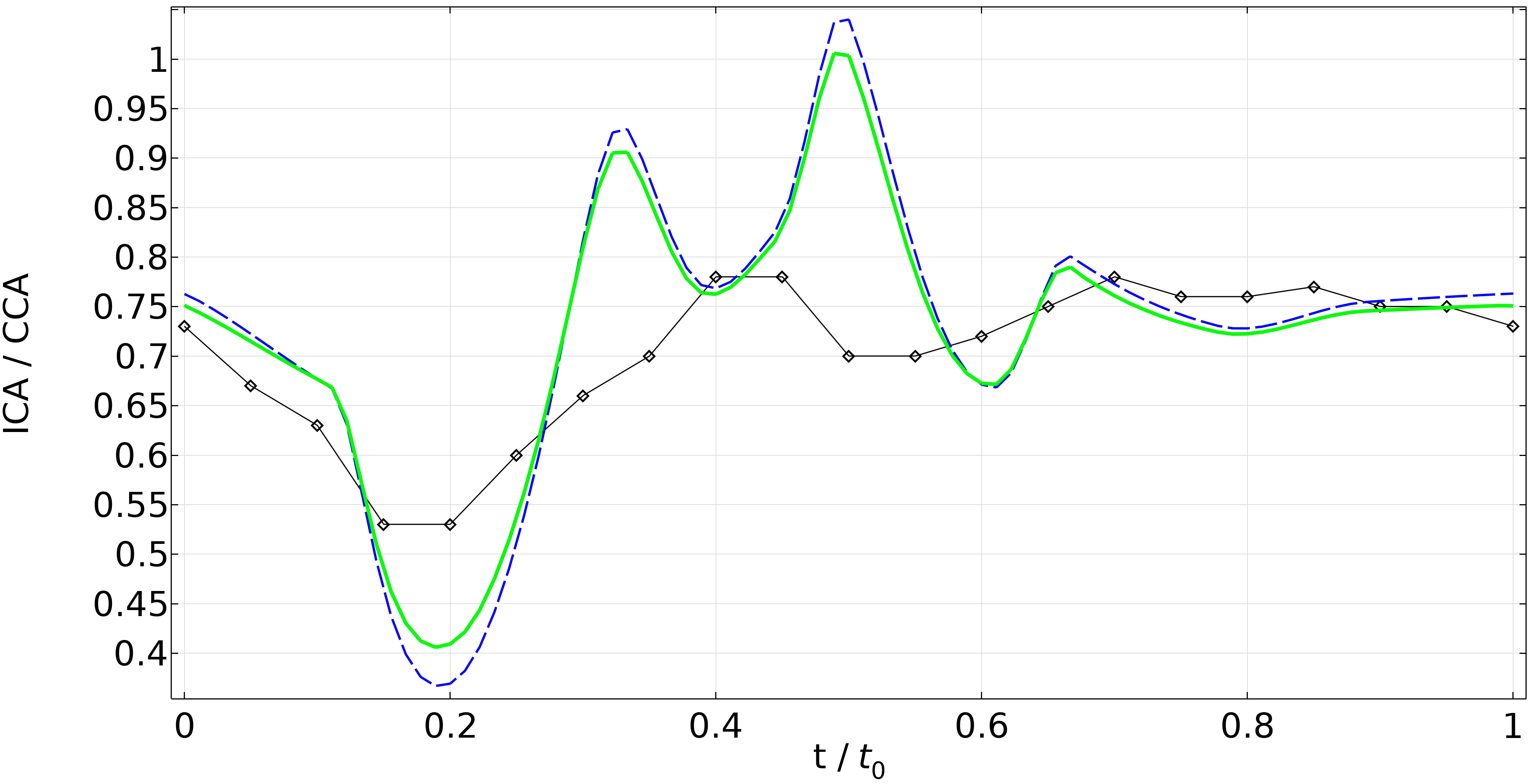}\\[1ex]
        \includegraphics[width=0.75\linewidth,height=5cm]{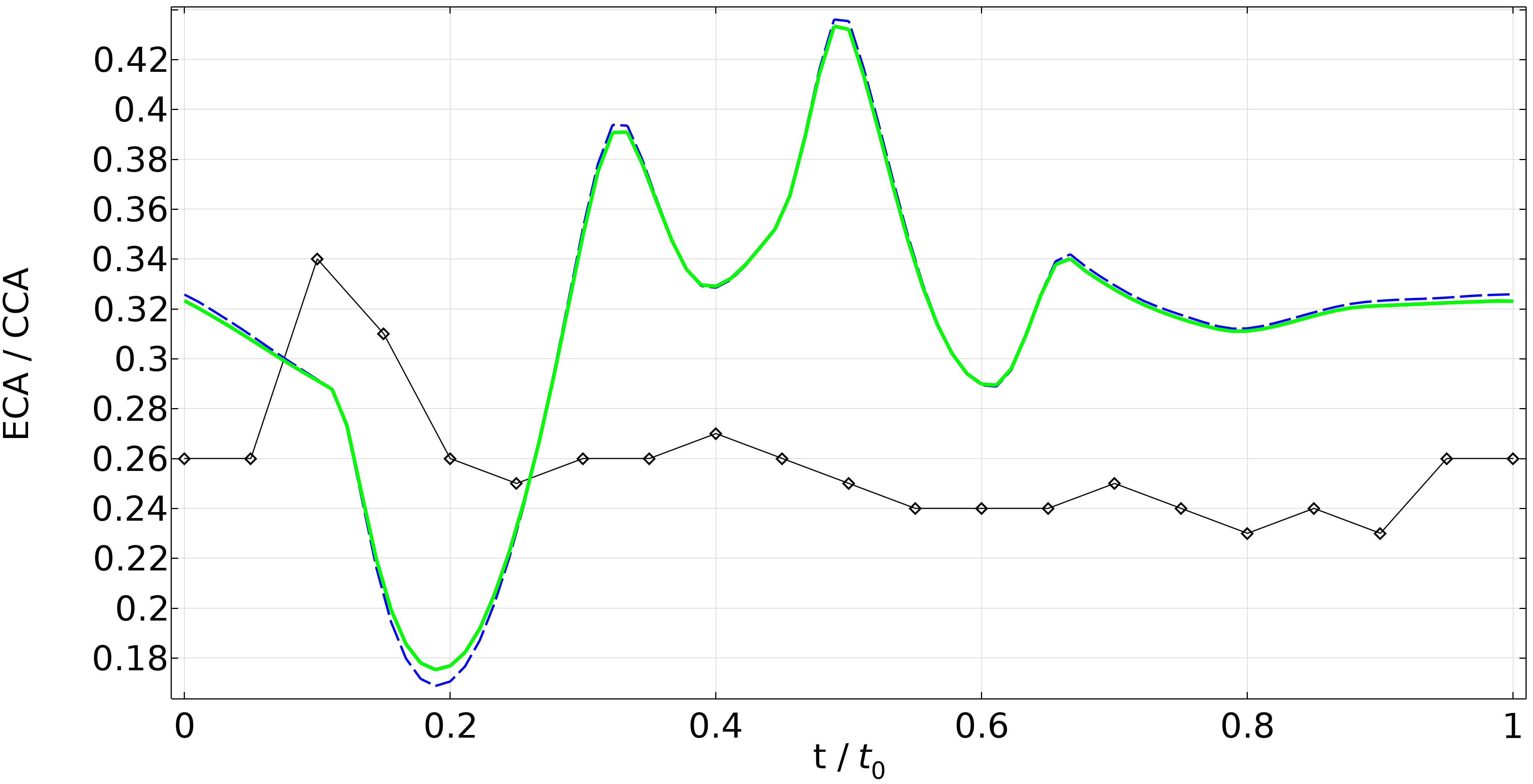}\\[1ex]
        \includegraphics[width=0.75\linewidth,height=5cm]{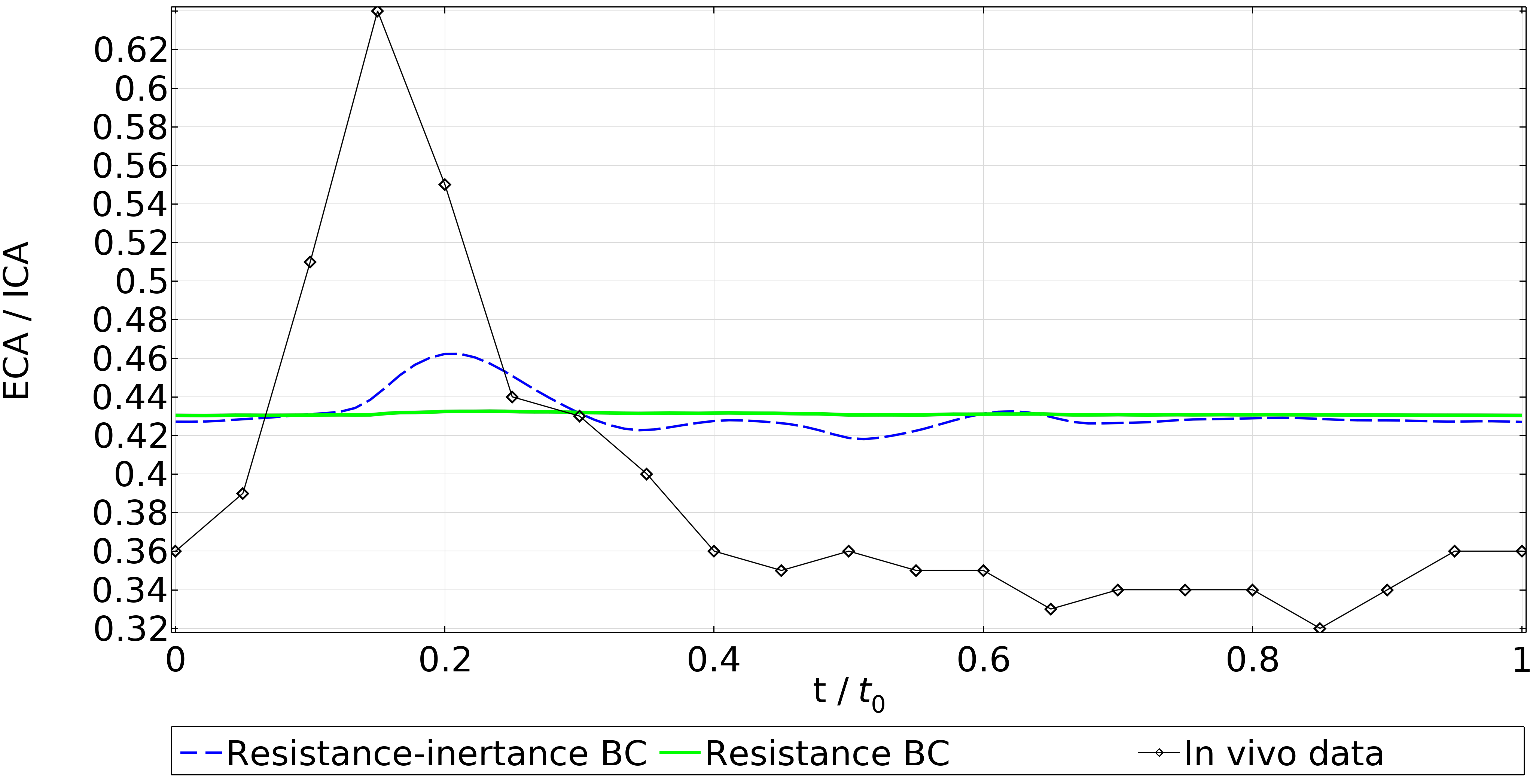}
    
        \caption{
         Temporal evolution of the flow ratios ICA/CCA (upper), ECA/CCA (middle), and ECA/ICA (lower) over one normalized cardiac cycle $t / t_0$, with $t_0=0.9$. Comparison of simulated ratios using resistance-inertance BC (\ref{eq:pout}) (blue dotted line) and pure resistance BC (green solid line) against clinical cohort data (black curves) from \cite[Fig. 4]{Marshall_2004}.
        }
        \vspace{0.05cm}
        \label{fig:flow_splitting}

\end{figure}
In Fig. \ref{fig:inflow} more detailed results for the time development of the pressure defined by \eqref{eq:pout} on the outlet surface is presented, divided into an resistance and an inertance part. The results demonstrate high resistance and very low inertial proportions of the boundary pressure, but with a significant systolic peak present in the inertial part of pressure for the ICA followed by slight temporal variation in diastolic phase. 
In Fig. \ref{fig:flow_splitting}, the resulting volumetric flow through the ICA and the collective flow through the ECA branches, normalized by the overall CCA inflow, are compared to clinical cohort data from \cite{Marshall_2004} as well as to simulation results using pure resistance boundary conditions (corresponding to $L_i=0$ in (\ref{eq:pout})) over one normalized cardiac cycle.
The course of the resulting simulated (instantaneous) ICA/CCA ratio (upper picture) shows similarities with clinical cohort data, but is clearly overestimated in the initial diastolic phase, which can be accounted to missing compliance effects of the downstream region which are not represented by the simplified inertance-resistance boundary condition (\ref{eq:pout}).
The ECA/ICA ratio (bottom graphic in Fig. \ref{fig:flow_splitting}) starts to capture the influence of inertance during systole when using the boundary condition (\ref{eq:pout}), unlike pure resistance boundary conditions which produce a constant flow ratio. However, our model of resistance-type BC (\ref{eq:pout}) still underestimates the strong increase of flow into the ECA compared to the noteworthy increase observed for the clinical ECA/ICA flow ratio in the systolic phase. The clinical ECA/CCA ratio (middle graphics in Fig. \ref{fig:flow_splitting}) is not represented by  resistance-inertance boundary condition (\ref{eq:pout}) at all, which is another indicator of the necessary improvement of the outlet boundary conditions e.g. by extending to Windkessel models.
One should also remark that variability in the morphology of ECA branches across the patient cohort may also has a non negligible influence on the observed discrepancy, especially on the ECA flow ratio for chosen individual CA morphology investigated in this study.
%

%%%%%%%%%%%%%%%%%%%%%%%%%%%%%%%%%%%%%%%%%%%%%%%%%%%%%%%%%%%
\subsection{Generalized Young's modulus: validation for clinical data }
\label{subsec:Invivo}
%%%%%%%%%%%%%%%%%%%%%%%%%%%%%%%%%%%%%%%%%%%%%%%%%%%%%%%%%%%
In the this section we demonstrate the feasibility of the above described modeling approach for strain-dependent Young's moduli from section \ref{subsec:young_moduli} to represent physiological conditions. 
The stress-strain behavior of silicon CA phantom, see Fig. \ref{fig:tensile_ds}, and the Young's modulus function derived from these  in vitro data, see Fig. \ref{fig:YoungsModulus}, does not fully replicate typical physiological behavior of arterial tissue, which is characterized by strain stiffening (collagen phase) following an initial elastin phase \cite{butlin2015age}. Instead, as shown in Fig. \ref{fig:YoungsModulus}, the curve originates in a stiffer state, exhibits strain softening within the physiological strain range, and demonstrates only marginal strain hardening at large strains. 
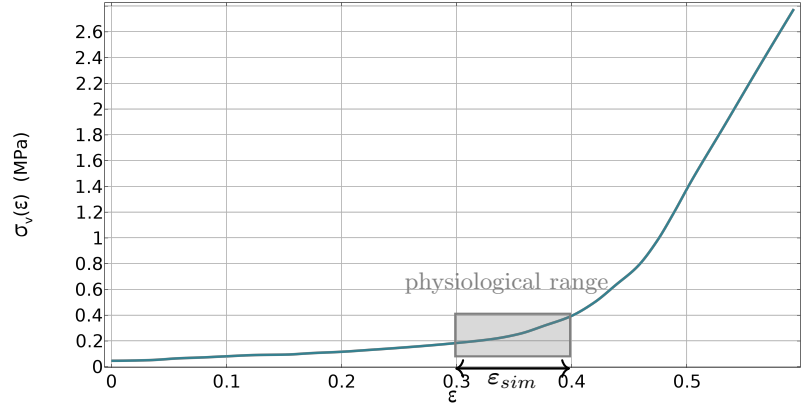
\begin{figure}[ht!]
    \captionsetup{justification=centering, margin=0.5cm}
    \centering
        \begin{tikzpicture}
            \node[anchor=south west, inner sep=0] (image) at (0,0) {\includegraphics[width=0.8\linewidth]{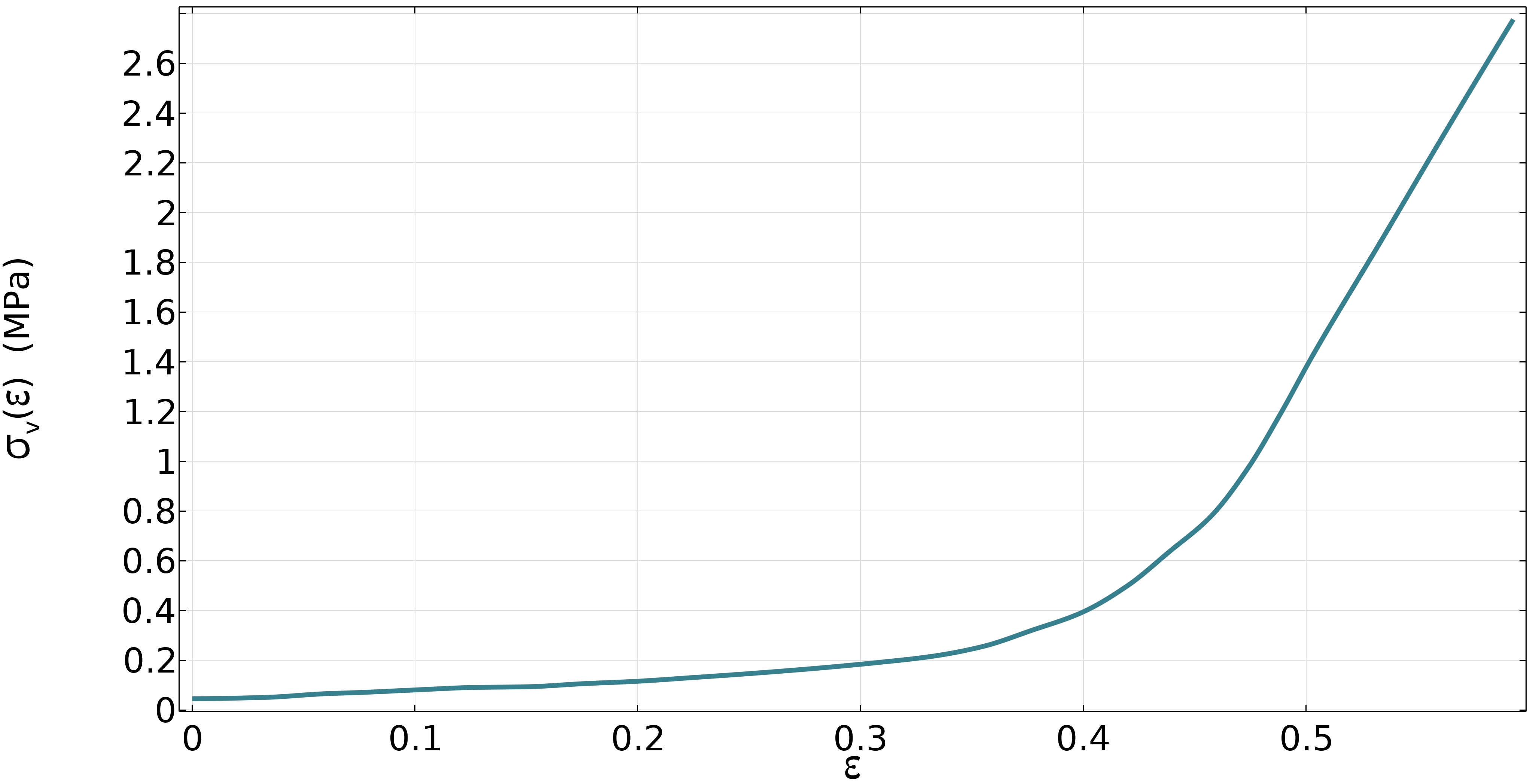}};
            \begin{scope}[x={(image.south east)},y={(image.north west)}]
                \draw[gray, thick, fill=gray, fill opacity=0.3] 
                    (0.56, 0.12) rectangle (0.705, 0.225);
                \node[gray, font=\scriptsize] at (0.625, 0.3) {physiological range};
                \draw[black,<->, thick] (0.56,0.09) -- (0.705,0.09);
                \node[black, font=\scriptsize] at (0.6325, 0.06) {$\boldsymbol{\varepsilon}_{sim}$};
            \end{scope}
        \end{tikzpicture}
        \caption{ Clinical data based stress-strain behavior of $\sigma_v (\boldsymbol{\varepsilon})$ for human CCA. Adapted from \cite{faturechi2019mechanical}, Fig. 1.\\}
        \vspace{0.4cm}
        \label{fig:hashemi_stressstrain}
    \end{figure}
To account for the inter-patient variability in stress-strain behavior \cite{boekhoven2015mechanical}  and to illustrate the potential of our Young's modulus approach derived from such stress-strain curves, we employ and validate a Young's modulus based on clinical stress-strain data for CCA from \cite{faturechi2019mechanical}.
For the realistic simulation of this physiological case we use the fact that arteries are commonly described as nearly incompressible materials - Yossef et al. \cite{yossef2017compressibility} reported a volumetric strain of approximately 5 \% in porcine common carotid arteries under physiological pressure loading $50- 200$ mmHg -  and bound the volumetric strain component of our strain metric $|\boldsymbol{\varepsilon}|$ (\ref{eq:em}) accordingly, s.t. $\frac{\varepsilon_{vol}}{3}\frac{1}{1-2\nu}\le0.05$ applies and is implemented in the strain metric $|\boldsymbol{\varepsilon}|$ of the generalized Young's modulus $E(|\boldsymbol{\varepsilon}|)$.

    \begin{figure}[ht!]
        \captionsetup{justification=centering, margin=0.5cm}
        \centering
        \includegraphics[width= 0.8\linewidth]{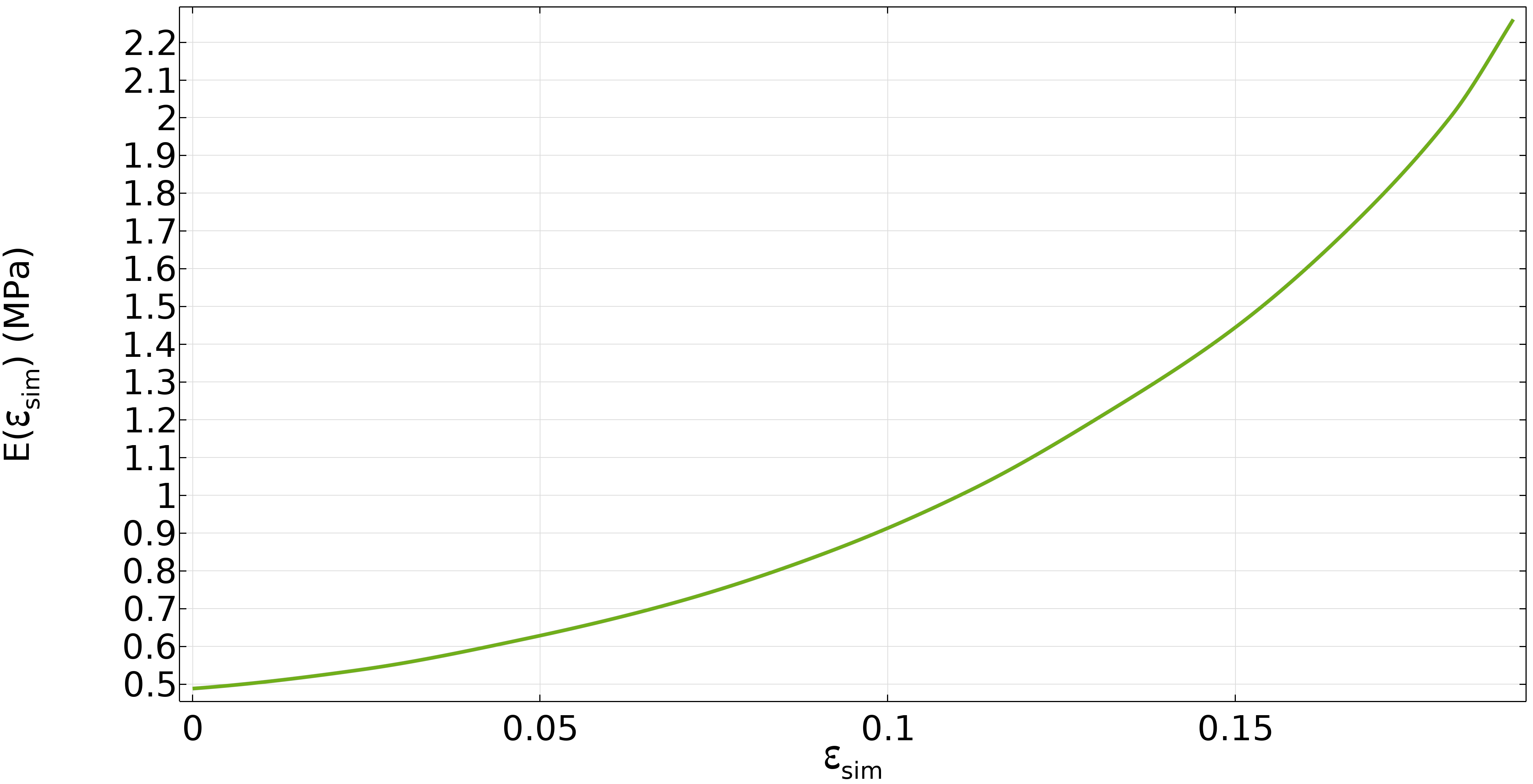}
        \caption{ Young's modulus function $E(\boldsymbol{\varepsilon}_{sim})$ extracted from tangential slopes of the clinical stress curve $\sigma_v (\boldsymbol{\varepsilon})$ from \cite{faturechi2019mechanical} for the physiological pressure range, see Fig. \ref{fig:hashemi_stressstrain}.}
        \label{fig:hashemi_youngs_physiological}
\end{figure}
The stress-strain curve reconstructed from clinical CCA data \cite{faturechi2019mechanical} is presented in Fig. \ref{fig:hashemi_stressstrain}
The curve displays a nearly linear response up to $30\%$ strain, after which it transitions to a strongly nonlinear regime. The physiological strain region for CCA is specified as $30\% - 40\%$ in \cite{faturechi2019mechanical}.
Since the patient-specific, from clinical imaging extracted CCA geometry already represents a prestrained/prestressed configuration under physiological conditions, the strain-dependent Young’s modulus for simulation was derived (similarly as in section \ref{subsec:young_moduli}) for strains $\varepsilon\ge30 \%$, assuming a spatial uniform prestrain of $0.3$ incorporated in the extracted geometry.
The Young's modulus curve for simulating the on clinical data based case, depicted in Fig. \ref{fig:hashemi_youngs_physiological}, therefore represents the modulus for the simulated strain $\boldsymbol{\varepsilon}_{sim}$, which contribute together with the prestrain to the total physiological strain $\boldsymbol{\varepsilon}=\boldsymbol{\varepsilon_{sim}}+0.3$. Note that the prestrain ($0.3$) is already included in the extracted CCA geometry and therefore in simulation only strains exceeding the prestrain range ($\boldsymbol{\varepsilon_{sim}}$) needs to be considered in the argument of $E(|\boldsymbol{\varepsilon}|)$. The resulting curve of $E$ exhibits a characteristic nonlinear behavior: Initially, elastin fibers dominate, resulting in a low tangential modulus, followed by progressive collagen fiber alignment, which induces strain-hardening at higher deformations.

To be consistent with the above assumption of a prestressed geometry extracted through clinical imaging, Bäumler et al. \cite{baeumler2020prestress} demonstrated the necessity of incorporating an initial prestress matrix $\mathbf{S}_0$ into the governing equations to prevent the overestimation of vessel deformation in fluid–structure interaction simulations. To achieve this, $\boldsymbol{S}_0$ has to be chosen in a way that it balances the fluid traction $\tilde{\boldsymbol{T}_f}$ exerted on the reference fluid-structure interface for the initial flow rate of $5 \ ml/s$ (see Fig. \ref{fig:inflow}, $t=0$), s.t.
\begin{equation} \label{eq:prestress}
    \tilde{\boldsymbol{T}}_f\boldsymbol{n}_f=-\boldsymbol{S}_0\boldsymbol{n}_s  \quad on \ \  \Gamma^0_{fsi}.
\end{equation}
Note that in this balanced state the deformation is zero and therefore $J=1$ and $\boldsymbol{F}=\boldsymbol{I}$ in \eqref{eq:prestress}, compare \eqref{eq:couplings}.
To account for this initial stress the solid deformation subproblem (\ref{eq:soliddeformation}), see \ref{subsec:sealedsetup}, becomes
\begin{equation}\label{eq:prestresssolid}
     \rho_s \frac{\partial^2 \boldsymbol{d}}{\partial t^2} = \nabla \cdot (\boldsymbol{P}+\boldsymbol{FS}_0)^T\quad\text{ in } \Omega^s_0 ,
\end{equation}
{and the fluid-structure interaction coupling condition, that balance the stresses of the fluid and the solid material (\ref{eq:couplings}) turn into the following boundary load}
\begin{equation}\label{eq:prestresscoupling}
 J\tilde{\boldsymbol{T}}_f\,\boldsymbol{n}_f= -(\boldsymbol{P}+\boldsymbol{FS}_0)^T \boldsymbol{n}_s \quad \text{on }\Gamma^0_{fsi}.
\end{equation}
To determine the prestress matrix $\boldsymbol{S}_0$ that satisfies condition~\eqref{eq:prestress}, we implement similar iterative procedure as described by Bäumler et al. \cite{baeumler2020prestress} in COMSOL Multiphysics. First, we approximate the fluid traction vector $\boldsymbol{h} := \tilde{\boldsymbol{T}}_f \boldsymbol{n}_f$ at the fluid-structure interface $\Gamma_{{fsi}}$ by running a rigid-wall stationary CFD simulation at the initial flow rate, solving the stationary Navier-Stokes equations~\eqref{eq:ns} without the time derivative term, decoupled from the solid subproblem. Next, the prestress tensor $\boldsymbol{S}_0$ is iteratively computed by solving the solid subproblem~\eqref{eq:prestresssolid} for $\boldsymbol{d}^{i+1}$ and $\boldsymbol{p}^{i+1}$ with boundary load $(\boldsymbol{P}^{(i+1)} + \boldsymbol{F}^{(i+1)}\boldsymbol{S}^{(i)}_0)^T \boldsymbol{n}_s = J\boldsymbol{h}$ on $\Gamma_{{fsi}}$. Starting with $\boldsymbol{S}^{(0)}_0 = \boldsymbol{0}$, the tensor is updated at each step via $\boldsymbol{S}^{(i+1)}_0 = (\boldsymbol{P}^{(i+1)} + \boldsymbol{F}^{(i+1)}\boldsymbol{S}^i_0)^T$ for iterations $i = 0, \dots, N$, here $N=40$. Finally, the resulting prestress tensor $\boldsymbol{S}_0 = \boldsymbol{S}^{(N)}_0$ is substituted into Eqs.~\eqref{eq:prestresssolid} and \eqref{eq:prestresscoupling} and the prestressed FSI problem \eqref{eq:ns}, \eqref{eq:prestresssolid}, \eqref{eq:couplings} with the prestress adaption of the dynamic coupling condition \eqref{eq:prestresscoupling} is computed.

\begin{figure}[htbp]
    \captionsetup{justification=centering, margin=1.5cm}
     \centering
     \begin{tikzpicture}
      \node[anchor=south west, inner sep=0] (img) at (0,0) {
            \adjincludegraphics[width=\linewidth]{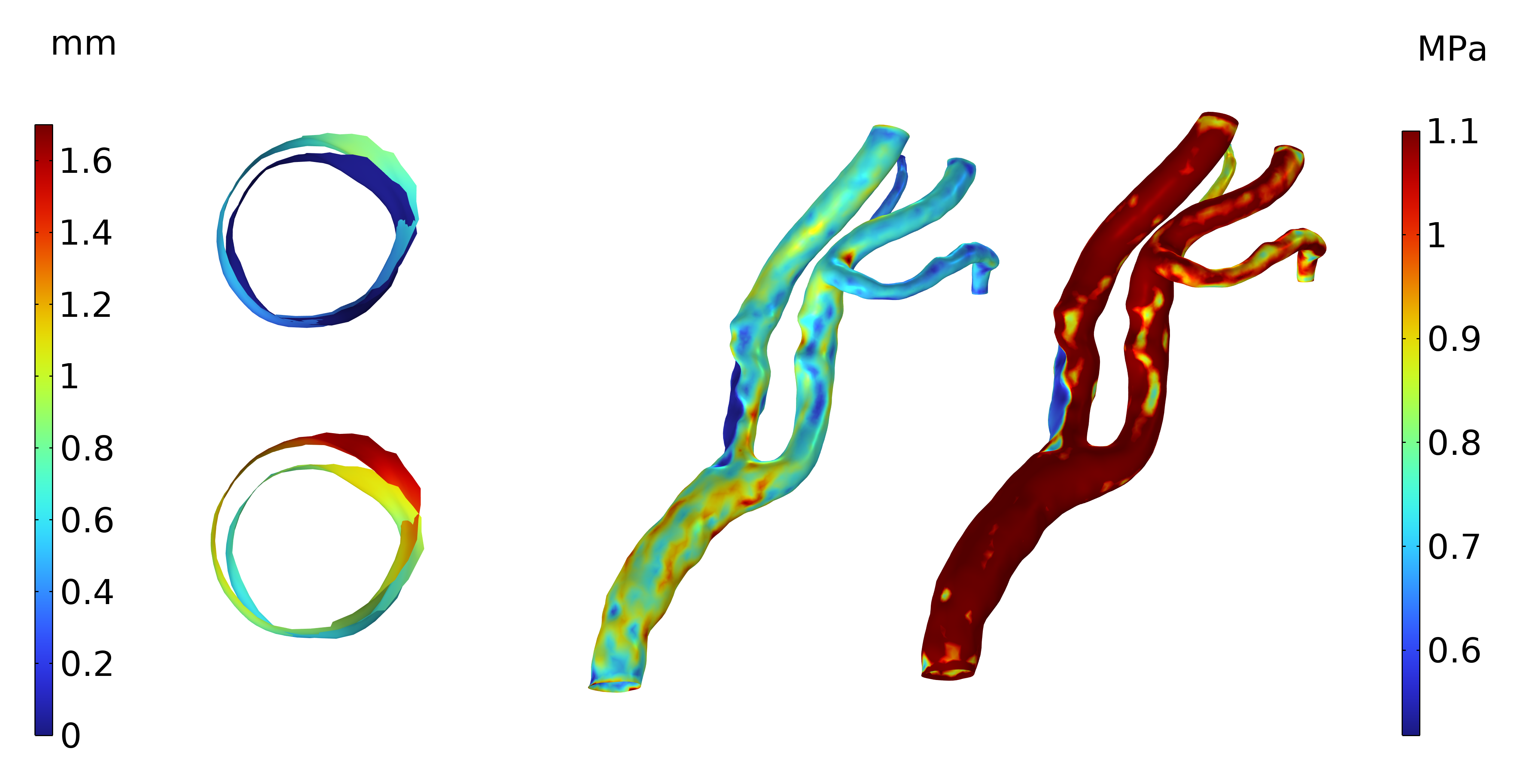}};
            \begin{scope}[x={(img.south east)}, y={(img.north west)}]
            % \draw[help lines, red, step=0.1] (0,0) grid (1,1); \foreach \x in {0,0.2,...,1} \node[red] at (\x,0) {\tiny \x}; \foreach \y in {0,0.2,...,1} \node[red] at (0,\y) {\tiny \y};
             \tikzset{
                myarrowt/.style={-{Stealth[scale=0.9]}, thick, draw=black, rounded corners=3pt},
                labelstyle/.style={ fill=white, fill opacity=0.85, text opacity=1, inner sep=2pt}
            }
            \tikzset{
                myarrowtt/.style={-{Stealth[scale=0.9]}, thick, draw=blue, rounded corners=3pt},
                labelstyle/.style={ fill=white, fill opacity=0.85, text opacity=1, inner sep=2pt}
            }
              \tikzset{
                myarrowttt/.style={-{Stealth[scale=0.9]}, thick, draw=gray, rounded corners=3pt},
                labelstyle/.style={ fill=white, fill opacity=0.85, text opacity=1, inner sep=2pt}
            }
            \node[labelstyle, align=center] (with) at (0.58, 0.95) {with \\ prestress};
            % \draw[myarrowt]  (with.east) -- (0.67,0.49);
            \node[labelstyle, align=center] (without) at (0.8, 0.95) {without \\ prestress};
            % \draw[myarrowt]  (without.east) -- (0.4,0.9);
            \node[labelstyle, align=center] (Diastole) at (0.2, 0.9) {Diastole};
            \node[labelstyle, align=center] (Systole) at (0.2, 0.5) {Systole};
            \node[labelstyle, align=center] (time) at (0.66, 0.89) {$t=1.1 \textbf{s}$};
            \node[labelstyle, align=center] (e_esp) at (0.52, 0.05) {$E(|\boldsymbol{\varepsilon}|)$};
            \node[labelstyle, align=center] (disp) at (0.2, 0.05) {Deformation ($\boldsymbol{d}$)};
            \node[labelstyle, align=center] (hash) at (0.37, 0.45) {without \\ prestress};
            \node[labelstyle, align=center] (star) at (0.37, 0.65) {with \\ prestress};
             \draw[myarrowt]  (star.west) --(0.15,0.7);
             \draw[myarrowt]  (star.west) --(0.15,0.3);
             \draw[myarrowtt]  (hash.west)--(0.26,0.78);
             \draw[myarrowtt]  (hash.west)--(0.25,0.42);
             \draw[rounded corners=3pt,gray, thick, fill=gray, fill opacity=0.2] 
                    (0.4, 0.3) rectangle (0.5, 0.35);
            \draw[gray, thick, fill=gray, fill opacity=0.05] 
                    (0.1, 0.1) rectangle (0.3, 0.95);
            \draw[myarrowttt]  (0.4,0.32) -- (0.3,0.32);
        \end{scope}
     \end{tikzpicture}
     \caption{Comparison of simulations initialized with and without prestress using $E(|\boldsymbol{\varepsilon}|)$ from Fig. \ref{fig:hashemi_youngs_physiological}.  \textbf{Left}: Deformation on a subregion of $\Gamma_{fsi}$ for diastolic and systolic pressure level. \textbf{Right}: Surface distribution   of $E(|\boldsymbol{\varepsilon}|)$ on $\Gamma_{fsi}$ at systole, $t=1.1s$.
     }
     \label{fig:e(eps)surf_prestrdisp}
 \end{figure}
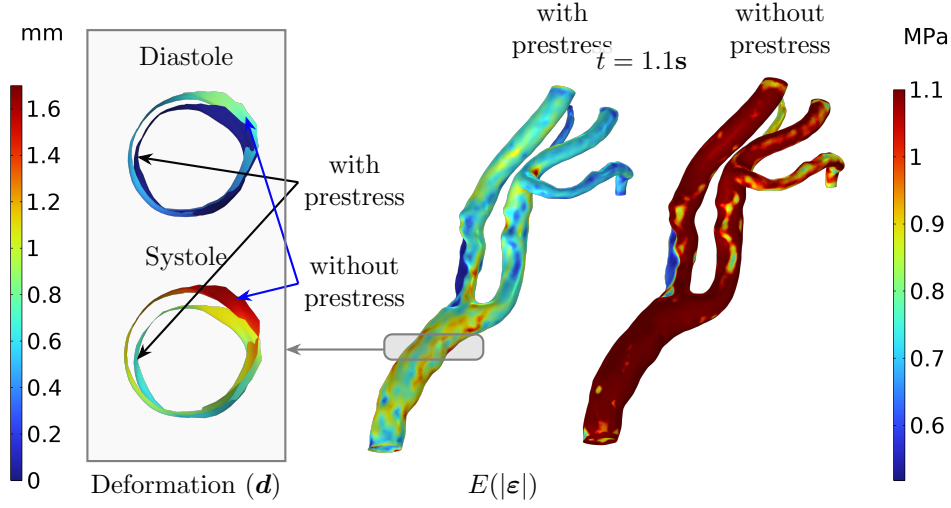

To evaluate the impact of prestress modeling on the proposed Young's modulus approach, the following simulation results incorporating prestress are also compared against simulations performed without prestress.
On the left of Fig. \ref{fig:e(eps)surf_prestrdisp} the deformation (on a subregion of $\Gamma_{fsi}$) of the simulations using the generalized Young's modulus based on clinical data, see Fig. \ref{fig:hashemi_youngs_physiological}, initialized with and without prestressing are compared at diastolic and systolic pressure level. For the simulation with prestress the deformation at diastolic pressure is zero, representing exactly the patient geometry extracted from CTA imaging, where as in the simulation without prestress deformations of $1 \text{ mm}$ are already present. The simulation with prestress is reaching a deformation of $1 \text{ mm}$ not until peak systole, where the deformation in the simulation without prestress  goes up to $1.6 \text{ mm}$ demonstrating the overestimation of the deformation when prestressing is not considered. 
The spatial distribution patterns of Young's modulus on $\Gamma_{fsi}$ at peak systole ($t=1.1\text{s}$) presented on the right of Fig. \ref{fig:e(eps)surf_prestrdisp} are similar for both cases using $E(|\boldsymbol{\varepsilon}|)$, but for the simulation without prestress the Young's modulus reach much higher amplitudes due to higher surface strains.
Interestingly, in both cases the stenotic plaque region experiences only minimal strains and therefore almost no increase in Young's modulus. This reduced strain capacity is attributed to the increased wall thickness in the stenotic region, which inherently limits deformation. Thus, no additional tuning of the Young’s modulus is  required to replicate the stiffer mechanical response of arterial plaque.
To assess the influence of strain-dependent elastic moduli on the simulation results, we compare
% the cases using in vitro (Figs.~\ref{fig:tensile_ds} and \ref{fig:YoungsModulus}) and
the generalized Young's modulus based on clinically observed stress–strain curve (Figs.~\ref{fig:hashemi_stressstrain} and \ref{fig:hashemi_youngs_physiological}) against a constant Young's modulus counterpart. For the simulation incorporating prestress and the generalized Young's modulus derived from clinical data the surface strain results in a systolic peak of approximately $0.1$, aligning well with the physiological range ($0.3-0.4$). Consequently, the equivalent constant Young's modulus is obtained by averaging $E(\boldsymbol{\varepsilon}_{sim})$ from Fig. \ref{fig:hashemi_youngs_physiological} over the strain interval $[0, 0.1]$, resulting in an Young's modulus of $0.64\text{ MPa}$.
We compare the temporal evolution and spatial distribution of the generalized Young's modulus, as well as volumetric distension and pressure–volume relationships.
\begin{figure}[htbp]
     \captionsetup{justification=centering, margin=0.5cm}
     \centering
      \begin{subfigure}[b]{0.95\textwidth}
    \hspace{0.3cm}
    \begin{tikzpicture}
        
        \node[anchor=south west,inner sep=0] (img) {\includegraphics[trim=0 0 0 0.2cm, clip,scale=1.3,width=\linewidth,height=1cm]{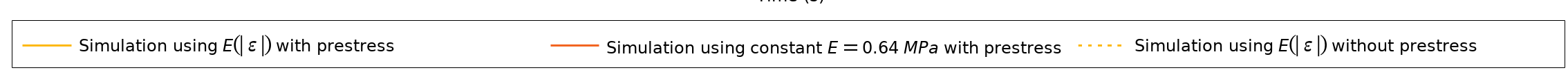}};

        \begin{scope}[x={(img.south east)}, y={(img.north west)}]
            
            \draw[white, fill=white] (0.05, 0.25) rectangle (0.3, 0.7);
            \draw[white, fill=white] (0.385, 0.25) rectangle (0.68, 0.7);
            \draw[white, fill=white] (0.72, 0.25) rectangle (0.95, 0.7);
            % 3. Neuer Text mit weißem Hintergrund
            \node[fill=white, inner sep=2pt,  scale=0.6,align=center] at (0.15, 0.45) {Simulation using $E(|\boldsymbol{\varepsilon}|)$ \\ with prestress};
             \node[fill=white, inner sep=2pt,  scale=0.6,align=center] at (0.55, 0.45) {Simulation using constant $E=0.64$ MPa \\ with prestress};
            \node[fill=white, inner sep=2pt,  scale=0.6,align=center] at (0.85, 0.45) {Simulation using $E(|\boldsymbol{\varepsilon}|)$ \\ without prestress};

        \end{scope}
    \end{tikzpicture}
    \end{subfigure}
    \begin{subfigure}[b]{0.49\textwidth}
        {\includegraphics[width=\linewidth,height=5cm]{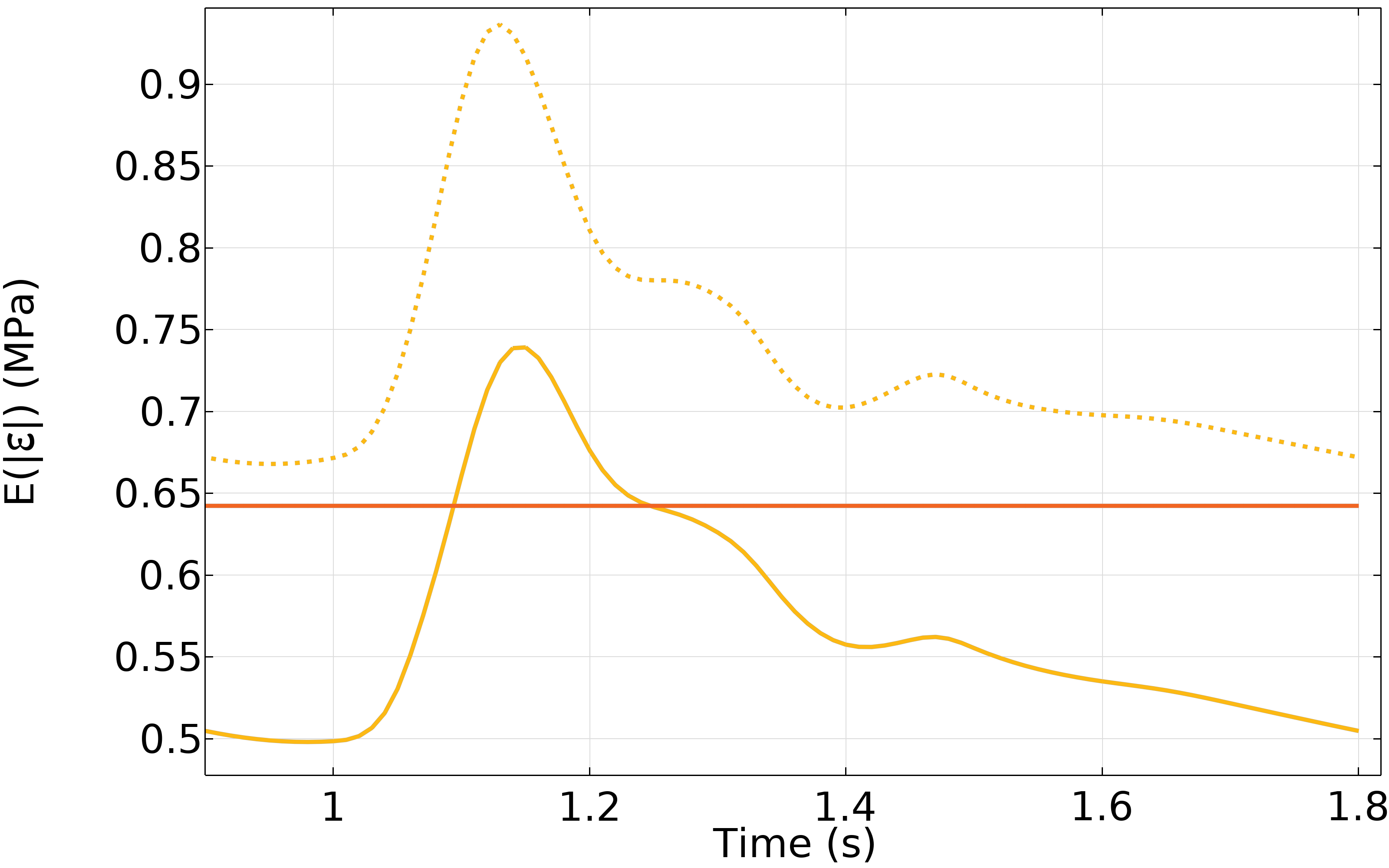}}
            \caption{  Time course of volume averaged generalized Young's moduli over the solid domain ($\Omega^s$).}
             \label{fig:e(eps)_comparison}
    \end{subfigure}
    \hfill
    \begin{subfigure}[b]{0.49\textwidth}
     \begin{tikzpicture}
          \node[anchor=south west, inner sep=0] (img) at (0,0) {%
            \adjincludegraphics[width=\linewidth,height=5cm]{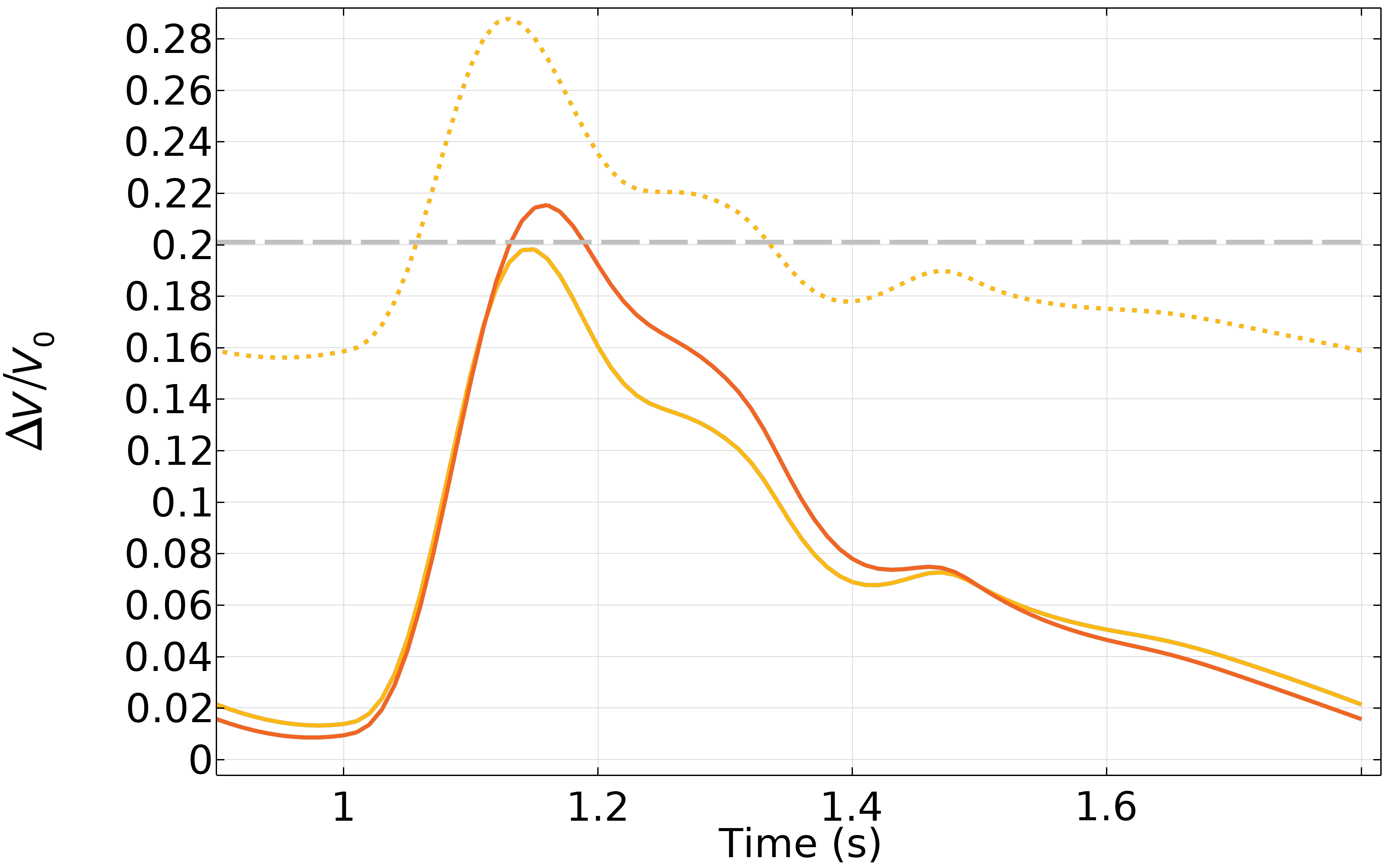}
            };
            \begin{scope}[x={(img.south east)}, y={(img.north west)}]

            \tikzset{
                myarrowt/.style={-{Stealth[scale=0.9]}, thin, draw=black, rounded corners=3pt},
                labelstyle/.style={ fill=white, fill opacity=0.85, text opacity=1, inner sep=2pt}
            }
            
            \node[gray, font=\scriptsize,align=center] at (0.65, 0.82) {systolic physiological \\ reference value: 0.201};
        \end{scope}
     \end{tikzpicture}
        \caption{  Relative volume change of the fluid domain $\Omega^f$ compared to systolic physiological reference value. }
        \label{fig:volume_growth}
    \end{subfigure}
    \caption{  Comparison of simulations initialized with prestress using constant $E$ and generalized $E(|\boldsymbol{\varepsilon}|)$ derived from clinical data, see Fig. \ref{fig:hashemi_youngs_physiological}, compared to simulation without prestress over one cardiac cycle.}
 \end{figure}
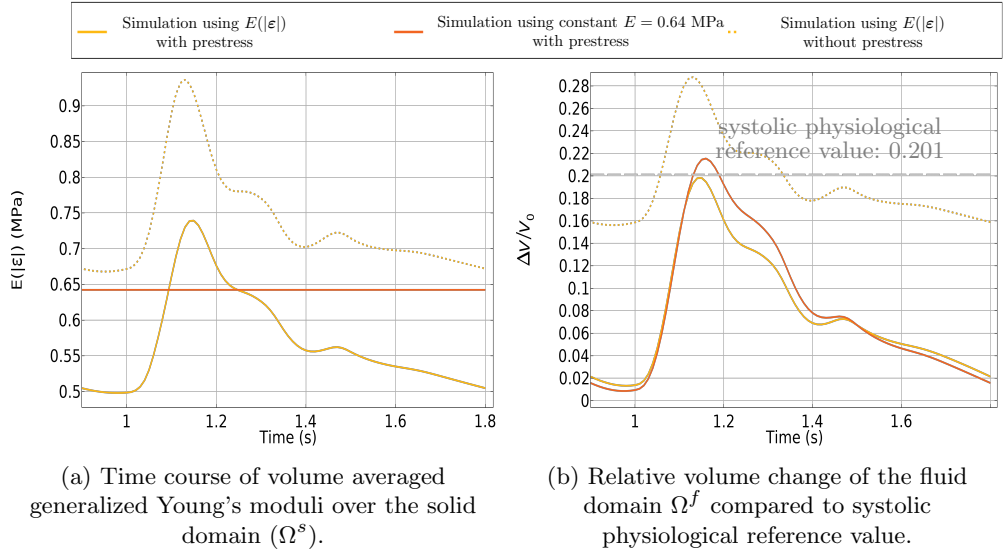
Fig. \ref{fig:e(eps)_comparison} illustrate the temporal evolution of the strain dependent Young’s modulus $E(|\boldsymbol{\varepsilon}|)$, averaged over the solid domain $\Omega^s$, for one cardiac cycle.
The simulation with prestressing demonstrates significantly Young's modulus  temporal variations ($0.5$ -- $0.74$ MPa), with pronounced stiffening during systole, reflecting the expected collagen fiber alignment in arterial tissue. The constant Young's modulus of $0.64$ MPa seems to be an appropriate average of this behavior but ignoring temporal and also
local spatial variations presented in Fig. \ref{fig:e(eps)surf_prestrdisp}. In the simulation without prestress the artery experience higher surface mean strains up to $0.13$ leading to much higher Young's modulus values due to the strongly non linear behavior for higher strains in Fig. \ref{fig:hashemi_youngs_physiological}.
\begin{wrapfigure}[19]{r}{0.3\linewidth}
    % \captionsetup{justification=centering, margin=1.5cm}
    \vspace{-0.2cm}
     \centering
 \begin{tikzpicture}
        \tikzset{
            myarrow/.style={-{Stealth[scale=1.2]}, thick, draw=black, rounded corners=3pt},
            labelstyle/.style={ fill=white, fill opacity=0.85, text opacity=1, inner sep=2pt}
        }
    \node[inner sep=0] (main) {\includegraphics[width=\linewidth]{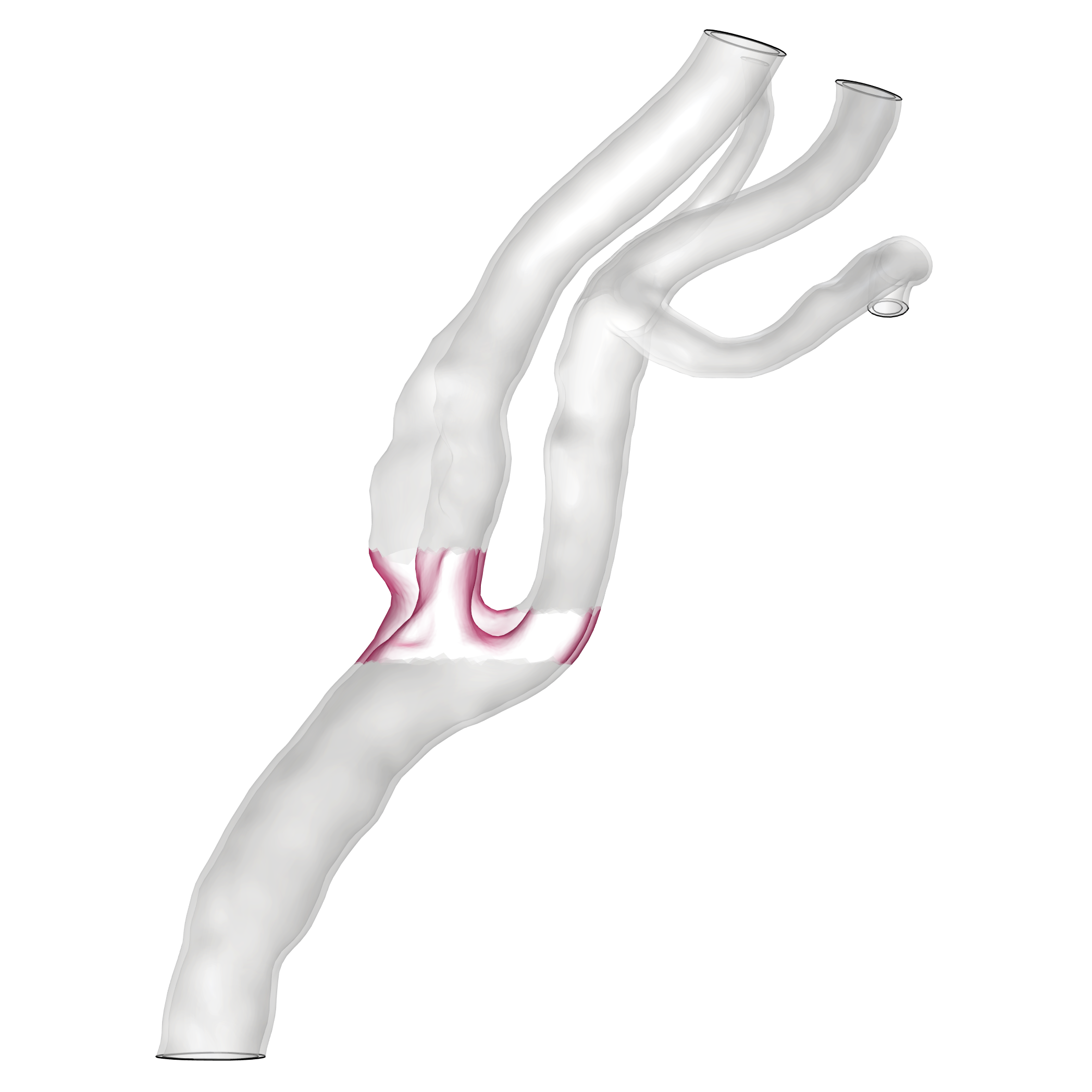}}; 
    \node[labelstyle, anchor=west] (ROI) at (-2,0.5) {ROI};
            \draw[myarrow] (-0.7,-0.45) --  (ROI)  ;
  \end{tikzpicture}
     \caption{  Region of interest (ROI) for the evaluation of simulated pressure-volume behavior matching the CCA sample region from \cite[Fig. 5.3]{boekhoven2015mechanical} for comparison.
     }
     \label{fig:ROI}
 \end{wrapfigure}
Fig. \ref{fig:volume_growth} displays the temporal development of the volumetric distension, $\Delta V/V_0 = (V - V_0)/V_0$, illustrating the relative volume expansion of the fluid domain $\Omega^{\text{f}}$ over one cardiac cycle for the three investigated cases. We compare the simulation results with a  physiological reference value of $(\Delta V/V_0)_{\text{phys}} \approx 20.1\%$ for peak systole, which is based on in vivo measured diameter distension values ($\Delta D/D_0 \approx 0.096$) for CCA from Segers et al. \cite{segers2004functional}. By assuming a constant arterial length over time (since the inlet and outlet positions are fixed in our model) and the cylindrical shape for CCA and each daugther branch, the volumetric distension corresponds to the relative cross-sectional area change, leading to: $(\Delta V/V_0)_{\text{phys}} \approx \Delta A/A_0 = 2(\Delta D/D_0) + (\Delta D/D_0)^2$. Including the above diameter distension value $0.096$ for all branches results in a value for $(\Delta V/V_0)_{\text{phys}}\approx 20.1\%$.
The maximum volumetric distension values for the prestressed simulations - $ \approx 0.215$ using a constant $E$ and $\approx 0.198$ using the generalized $E(|\boldsymbol{\varepsilon}|)$ - are in almost perfect agreement with the physiological reference value $0.201$. In contrast, simulation not including prestress results in a maximum value of $\approx 0.288$, which noticeably exceeds the physiological reference.

\begin{figure}[htbp]
    \captionsetup{justification=centering, margin=1.5cm}
     \centering
     \includegraphics[width=0.8\linewidth]{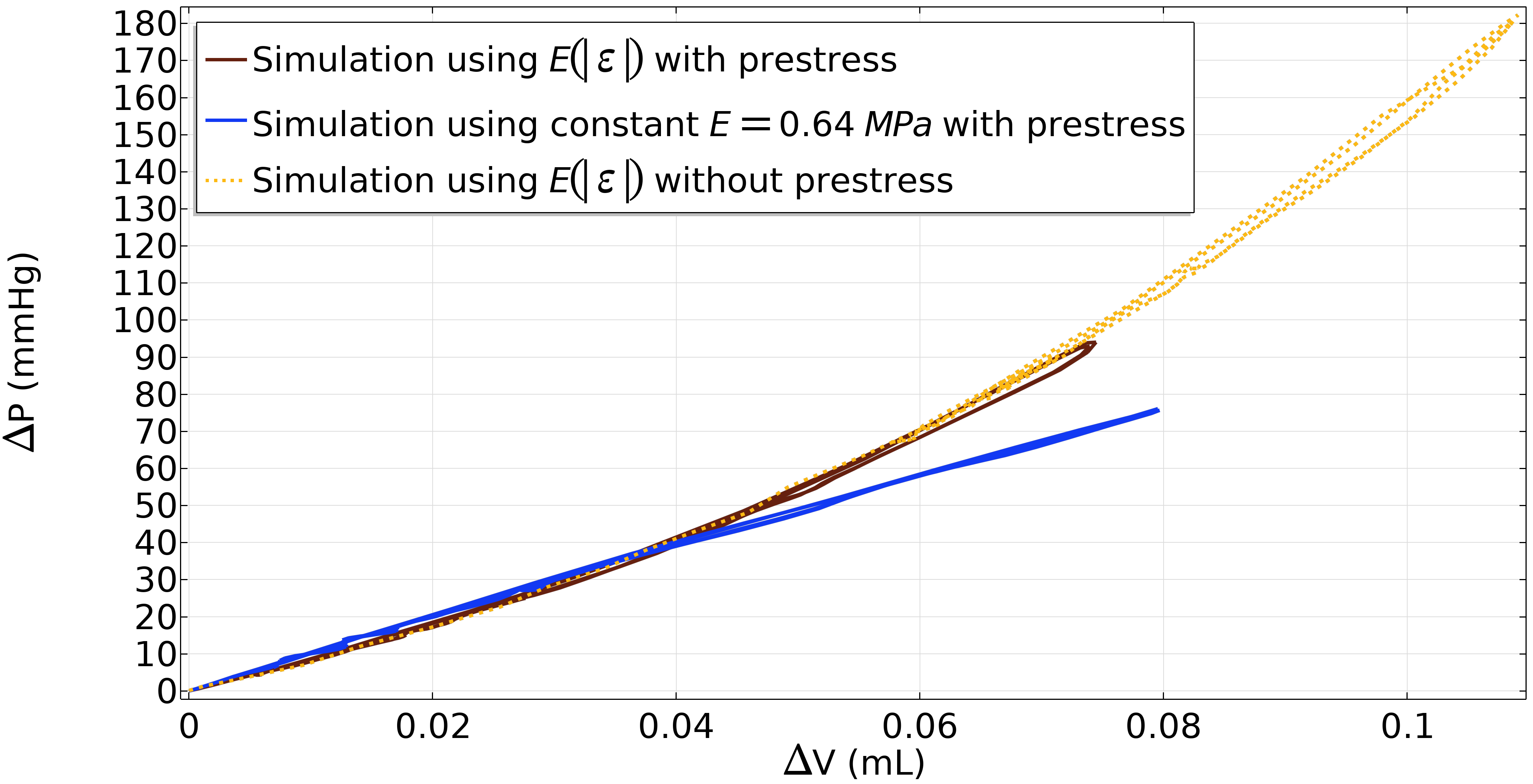}
     \caption{  Pressure growth to volume growth relationships for the ROI comparing constant and generalized linear elastic material models using $E(|\boldsymbol{\varepsilon}|)$.
     }
     \label{fig:realistic_pressure_volume}
 \end{figure}

Finally, Fig.~\ref{fig:realistic_pressure_volume} presents the pressure-volume relationships in terms of pressure change ($\Delta P = P - P_0$) and volume change ($\Delta V = V - V_0$) for the three considered cases. These were evaluated in a specific region of interest as a subregion of $\Omega^{f}$ ,see Fig.~\ref{fig:ROI}, enabling a direct comparison with published in vivo pressure-volume loops for CCA bifurcation samples from Boekhoven~\cite[Fig.~5.3]{boekhoven2015mechanical}.
  A clear improvement of the strain dependent Young's modulus approach compared to the constant one can be observed in Fig. \ref{fig:realistic_pressure_volume}. Where the latter show a pure linear pressure-volume response, the strain-dependency of $E$ lead to a pronounced nonlinear response aligning well with the measured behavior by Boekhoven \cite[Fig. 5.5]{boekhoven2015mechanical}. Here, for the generalized $E$ simulation with prestress we observe $0.07 \ \text{mL}$ volume increase induced through a pressure action of $90 \ \text{mmHg}$ matching the results measured by Boekhoven \cite[Fig. 5.5, Sample IX]{boekhoven2015mechanical}.  The simulation without prestressing show a similar pressure-volume response behavior but do reach (significantly) higher pressure and volume change values, compared to reported values.
  Noteworthy, for all cases high up to unphysiological systolic pressure levels - $140 \ \text{mmHg}$ using constant $E$, $160 \ \text{mmHg}$ using generalized $E$ with and $180 \ \text{mmHg}$ when no prestress applied - are present (for the sake of the papers extend these simulation results are not presented). This observed pressure overshoot can be attributed to the patient-unspecific resistance-inertance boundary condition, which neglects the pressure-damping effect of downstream vascular compliance. To represent the physiological conditions in more realistic clinical application a transition to more sophisticated Windkessel models is required, including patient-specific optimization of the compliance parameters, see eq. \eqref{eq:wk4}.
Nevertheless, our results illustrate clearly the capacity of our strain-dependent Young’s modulus approach to effectively captures diverse mechanical behaviors of arterial walls, thereby enhancing the physiological relevance of the simulations  without adapting the underlying elasticity model. The presented results also underline the importance of the incorporation of prestress modeling to prevent the overestimation of deformation and of proposed strain-dependent Young's modulus.
%

%-----------------------------------------------------------------------------------------------------------------------------------------------------------------------------------------------------------------------------------------------------
\section{Summary}
In this contribution, we present a numerical framework for modeling and validating the compliance behavior of an individual carotid artery tree using a  FSI approach.
We extend the conventional linear elasticity assumption by incorporating a strain-dependent Young’s modulus, derived from the underlying stress-strain data. The in vitro validation of volumetric compliance and diametric distensibility in simulations of a sealed carotid artery model demonstrates that the adapted elasticity formulation yields numerical results in closer agreement with the physical model than a constant Young’s modulus.
To simulate realistic cropped arterial geometries with outlet cross-sections and an unresolved downstream region in the absence of velocity data, we propose a resistance-inertance boundary condition. This approach, derived from more complex Windkessel models, accounts for the resistance and inertia of the downstream circulatory system within the FSI framework, but it disregards the impact of the downstream compliance.
In the cropped (open) artery the strain-dependent Young’s modulus model is validated against both in vitro measured data (as described in section \ref{subsec:laboratorysetup}) as well as for clinically measured stress-strain data for a common carotid artery (CCA) from \cite{faturechi2019mechanical}, using a constant Young’s modulus as a reference. The comparison confirms that our approach accurately reproduces distinct mechanical wall behaviors, offering a novel method to support and validate the development of in vitro arterial models (e.g., thin wall as well as artery block models \cite{shiravand2026,shiravand2025core}). Furthermore, it enables the replication of patient-specific wall mechanics within a linear elasticity framework without the necessity of non-linear elasticity solvers. The physiological relevance of the presented numerical study is further underscored by incorporating prestress into the model reflecting more accurately the physiological state in terms of strain, volumetric distension and pressure-volume response and finally vessel compliance and deformation then in simulations without prestress. 

A current limitation of this study is the presented outlet resistance-inertance boundary condition, which neglects the compliance effects of the downstream vasculature. Whereas this outlet conditions are appropriate to model the global compliance behavior of the vessel more precisely w.r. to measurements, than e.g., zero stress free outflow condition, they  result in unphysiologically high systolic pressure peaks - particularly in cases exhibiting strong nonlinear strain-hardening behavior - and in non physiological flow split into the daughter branches.
Therefore, the future work should focus on enhancing the proposed numerical model and investigation of more sophisticated Windkessel boundary conditions (\ref{eq:wk4}). This involves adaptation and optimization of existing Windkessel models for patient specific modeling within the FSI framework.

%-----------------------------------------------------------------------------------------------------------------------------------------------------------------------------------------------------------------------------------------------------

\backmatter
% Please refer to Journal-level guidance for any specific requirements.
\bmhead{Acknowledgements}
This work was partially funded by the BMBF joint Project \href{https://math4innovation.de/index.php?id=70}{05M20UNA-MLgSA}.
% \section*{Declarations}

% Some journals require declarations to be submitted in a standardised format. Please check the Instructions for Authors of the journal to which you are submitting to see if you need to complete this section. If yes, your manuscript must contain the following sections under the heading `Declarations':

% \begin{itemize}
% \item Conflict of interest/Competing interests (check journal-specific guidelines for which heading to use)
% \item Ethics approval and consent to participate
% \item Consent for publication
% \item Materials availability
% \end{itemize}


\bibliography{sn-bibliography}% common bib file

\end{document}